\documentclass[12pt,openany]{memo-l}
\usepackage{amssymb, amstext, amscd, amsmath}

\theoremstyle{plain}
\newtheorem{thm}{Theorem}[chapter]
\newtheorem{cor}[thm]{Corollary}
\newtheorem{prop}[thm]{Proposition}
\newtheorem{lem}[thm]{Lemma}

\theoremstyle{definition}
\newtheorem{rem}[thm]{Remark}

\newtheorem{defn}[thm]{Definition}

\newtheorem{eg}[thm]{Example}
\newtheorem{conj}[thm]{Conjecture}

\numberwithin{section}{chapter}
\numberwithin{equation}{chapter}

\newcommand{\bB}{{\mathbb{B}}}
\newcommand{\bC}{{\mathbb{C}}}
\newcommand{\bD}{{\mathbb{D}}}
\newcommand{\bF}{{\mathbb{F}}}
\newcommand{\bN}{{\mathbb{N}}}

\newcommand{\bR}{{\mathbb{R}}}
\newcommand{\bT}{{\mathbb{T}}}
\newcommand{\bZ}{{\mathbb{Z}}}

  \newcommand{\A}{{\mathcal{A}}}
  \newcommand{\B}{{\mathcal{B}}}
  \newcommand{\C}{{\mathcal{C}}}
  
  \newcommand{\E}{{\mathcal{E}}}
  \newcommand{\F}{{\mathcal{F}}}
  \newcommand{\G}{{\mathcal{G}}}
\renewcommand{\H}{{\mathcal{H}}}
  \newcommand{\I}{{\mathcal{I}}}
  \newcommand{\J}{{\mathcal{J}}}
  \newcommand{\K}{{\mathcal{K}}}
\renewcommand{\L}{{\mathcal{L}}}
  \newcommand{\M}{{\mathcal{M}}}
  \newcommand{\N}{{\mathcal{N}}}
\renewcommand{\O}{{\mathcal{O}}}
\renewcommand{\P}{{\mathcal{P}}}
  \newcommand{\Q}{{\mathcal{Q}}}
  
\renewcommand{\S}{{\mathcal{S}}}
  \newcommand{\T}{{\mathcal{T}}}
  \newcommand{\U}{{\mathcal{U}}}
  \newcommand{\V}{{\mathcal{V}}}
  
  \newcommand{\X}{{\mathcal{X}}}
  
  \newcommand{\Z}{{\mathcal{Z}}}

\newcommand{\ep}{\varepsilon}
\renewcommand{\phi}{\varphi}
\newcommand{\upchi}{{\raise.35ex\hbox{\ensuremath{\chi}}}}

\newcommand{\fA}{{\mathfrak{A}}}
\newcommand{\fD}{{\mathfrak{D}}}
\newcommand{\fK}{{\mathfrak{K}}}
\newcommand{\fL}{{\mathfrak{L}}}
\newcommand{\fM}{{\mathfrak{M}}}
\newcommand{\fT}{{\mathfrak{T}}}

\newcommand{\fs}{{\mathfrak{s}}}
\newcommand{\ft}{{\mathfrak{t}}}

\newcommand{\ba}{{\mathbf{a}}}

\newcommand{\bi}{{\mathbf{i}}}

\newcommand{\bk}{{\mathbf{k}}}

\newcommand{\rC}{{\mathrm{C}}}

\newcommand{\qand}{\quad\text{and}\quad}
\newcommand{\qif}{\quad\text{if}\quad}
\newcommand{\qfor}{\quad\text{for}\quad}
\newcommand{\qforal}{\quad\text{for all}\quad}

\newcommand{\AND}{\text{ and }}
\newcommand{\FOR}{\text{ for }}
\newcommand{\FORAL}{\text{ for all }}

\newcommand{\OR}{\text{ or }}

\newcommand{\Aut}{\operatorname{Aut}}
\newcommand{\diag}{\operatorname{diag}}
\newcommand{\End}{\operatorname{End}}
\newcommand{\id}{{\operatorname{id}}}
\newcommand{\Lat}{\operatorname{Lat}}

\newcommand{\spn}{\operatorname{span}}
\newcommand{\rad}{\operatorname{rad}}
\newcommand{\rep}{\operatorname{rep}}
\newcommand{\spr}{\operatorname{spr}}
\newcommand{\Tr}{\operatorname{Tr}}

\newcommand{\bsl}{\setminus}
\newcommand{\ca}{\mathrm{C}^*}
\newcommand{\cenv}{\mathrm{C}^*_{\text{env}}}
\newcommand{\dint}{\displaystyle\int}
\newcommand{\Fn}{\mathbb{F}_n^+}
\newcommand{\Ftwo}{\mathbb{F}_2^+}
\newcommand{\Fock}{\ell^2(\Fn)}
\newcommand{\lip}{\langle}
\newcommand{\rip}{\rangle}
\newcommand{\ip}[1]{\langle #1 \rangle}
\newcommand{\mt}{\varnothing}
\newcommand{\ol}{\overline}

\newcommand{\wot}{\textsc{wot}}

\newcommand{\ltwo}{\ell^2}
\newcommand{\sotsum}{\textsc{sot--}\!\!\sum}
\newenvironment{sbmatrix}{\left[\begin{smallmatrix}}{\end{smallmatrix}\right]}

\begin{document}
\frontmatter

\title[Operator algebras]{Operator algebras\\
for multivariable dynamics}

\author[K.R. Davidson]{Kenneth R. Davidson}
\address{Pure Mathematics Department\\University of \mbox{Waterloo}\\Waterloo, ON\;
N2L--3G1\\CANADA}
\email{krdavids@uwaterloo.ca}

\author[E.G. Katsoulis]{Elias~G.~Katsoulis}
\address{Department of Mathematics\\ East Carolina University\\
Greenville, NC 27858\\USA}
\email{KatsoulisE@mail.ecu.edu}

\thanks{First author partially supported by an NSERC grant.}
\thanks{Second author was partially supported by a grant from ECU}

\begin{abstract}
Let $X$ be a locally compact Hausdorff space with
$n$  proper continuous self maps $\sigma_i:X \to X$ for $1 \le i \le n$.
To this we associate two conjugacy operator algebras which
emerge as the natural candidates for the universal algebra
of the system, the tensor algebra $\A( X , \tau)$ and the semicrossed
product $\rC_0(X)\times_\tau\Fn$.

We develop the necessary dilation theory for both models.
In particular, we exhibit an explicit family of boundary representations
which determine the C*-envelope of the tensor algebra.

We introduce a new concept of conjugacy for multidimensional systems,
called piecewise conjugacy. We prove that
the piecewise conjugacy class of the system can be
recovered from the algebraic structure of either
$\mathcal A( X , \sigma)$ or $\mathrm{C}_0(X)\times_\sigma \mathbb{F}_n^+$.
Various classification results follow as a consequence.
For example, if $n=2$ or $3$, or the space $X$ has covering dimension
at most 1, then the tensor algebras are algebraically isomorphic
(or completely isometrically isomorphic) if and only if the systems are
piecewise topologically conjugate.

We define a generalized notion of wandering sets and recurrence.
Using this, it is shown that $\mathcal A( X , \sigma)$ or
$\mathrm{C}_0(X)\times_\sigma \mathbb{F}_n^+$
is semisimple if and only if there are no generalized wandering sets.
In the metrizable case, this is equivalent to each $\sigma_i$ being
surjective and $v$-recurrent points being dense for each $v \in \mathbb{F}_n^+$.
\end{abstract}

\date{}
\maketitle

\setcounter{page}{4}

\tableofcontents

\mainmatter
\chapter{Introduction}\label{C:intro}

Let $X$ be a locally compact Hausdorff space; and suppose we are given
$n$  proper continuous self maps $\sigma_i:X \to X$ for $1 \le i \le n$, i.e., a
\textit{multivariable dynamical system}.
In this paper, we develop a theory of conjugacy algebras for such
multivariable dynamical systems.
One of the goals is to develop connections between the dynamics of
multivariable systems and fundamental concepts in operator algebra theory.
One of the main outcomes of this work is that the classification and
representation theory of conjugacy algebras
is intimately connected to piecewise conjugacy and generalized recurrence
for multivariable systems.

In the case of a dynamical system with a single map $\sigma$,
there is one natural prototypical operator algebra associated to it,
the semicrossed product of the system.
As we shall see, the case $n >1$ offers a far greater diversity of examples.
It happens that there are various non-isomorphic algebras that can
serve as a prototype for the conjugacy algebra of the system.
The algebras should contain an isometric copy of $\rC_0(X)$;
plus they need to contain generators $\fs_i$ which encode the covariance
relations of the maps $\sigma_i$.
In addition, it is necessary to impose norm conditions to be able to
talk about a \textit{universal operator algebra} for the system.
The choice of these conditions creates two natural choices
for the appropriate universal operator algebra for the system:
the case in which the generators
are either isometric, producing the \textit{semicrossed product}, or
row isometric, producing the \textit{tensor algebra}.

The main goal of the paper is to demonstrate that these operator algebras
encode (most of) the dynamical system.
The strongest possible information that might be recovered from
an operator algebra of the form we propose would be to obtain the
system up to conjugacy and permutation of the maps, since there
is no intrinsic order on the generators.
It turns out that what naturally occurs is a local conjugacy, in which
the permutation may change from one place to another.
This leads us to the notion that we call \textit{piecewise conjugacy}.
We will show that either of our universal operator algebras contains
enough information to recover the dynamical system up to piecewise
conjugacy.

These results offer new insights into the classification theory for
operator algebras. In \cite{MS3}, Muhly and Solel initiated an
ambitious program of classifying all tensor algebras of
C*-correspondences up to isomorphism.
They introduced a notion of aperiodicity for C*-correspondences,
and were able to classify up to isometric isomorphism
all tensor algebras associated with aperiodic correspondences.
Many important operator algebras, including various natural subalgebras
of the Cuntz algebras, were left out of their remarkable classification scheme.
A first effort to address the periodic case was the study of isomorphisms
between graph algebras \cite{KK, KP, Sol}.
The classification results of this paper for tensor algebras of multidimensional
systems includes many examples which are not aperiodic,
and also pushes the envelope beyond isometric isomorphisms.
The complexity of the arguments involved in our setting,
as well as the need for importing non-trivial results from
other fields of mathematics, seem to indicate that a comprehensive
treatment of the periodic case for arbitrary tensor algebras of
C*-correspondences may not be feasible at this time.

As a first step in understanding our operator algebras for
multi-variable dynamical systems, we produce
concrete models by way of dilation theory.
In recent work, Dritschel and McCullough \cite{DMc} show that
the \textit{maximal representations} of an operator algebra $\A$ are precisely those
which extend (uniquely) to a $*$-representation of the C*-envelope, $\cenv(\A)$.
They use this to provide a new proof of the existence of the C*-envelope
independent of Hamana's theory \cite{Ham} of injective envelopes.
When such representations are irreducible, they are called boundary
representations.  This key notion was introduced by Arveson in
his seminal work \cite{Arv1} on dilation theory for operator algebras
(which are generally neither abelian nor self-adjoint).
Very recently, Arveson \cite{Arv_choq} has shown that there are always
sufficiently many boundary representations to determine $\cenv(\A)$.

The dilation theory is both more straightforward and more satisfying in the
case of the tensor algebra.
One can explicitly exhibit a natural and tractable family of
boundary representations which yield a completely isometric
representation of the operator algebra.
This provides the first view of the C*-envelope.

It also turns out that the tensor algebra is a C*-correspondence
in the sense of Muhly and Solel \cite{MS2}.
This enables us to exploit their work, and work of Katsura \cite{Ka}
and Katsoulis--Kribs \cite{KK3}, in order to describe the
C*-envelope of the tensor algebra as a Cuntz--Pimsner algebra.

In the semicrossed product situation, one needs to work harder to
achieve what we call a \textit{full dilation}.
These are the maximal dilations in this context.
This allows us to show that generally these algebras
are not C*-correspondences.
We have no `nice' class of representations which yield a completely
isometric representation.  So the explicit form of the C*-envelope
remains somewhat obscure in this case. Nevertheless, the information gained from
the dilation theory for the semicrossed product plays an important role in the sequel.
Indeed, in Example \ref{Ex:pwc but not iso} we use finiteness for the
$\ca$-envelope to show that the classification scheme for tensor algebras (see below) is
not applicable in the semicrossed product situation.

We then turn to the problem of recovering the dynamics from the
operator algebra.  As a first step, we establish that algebraic isomorphisms
between two algebras of this type are automatically continuous.
Then we apply the techniques from our analysis \cite{DKisotop}
of the one-variable case to study the space of characters and
the two-dimensional nest representations.
These spaces carry a natural analytic structure which is critical
to the analysis.
The fact that we are working in several variables means that we
need to rely on some well-understood but non-trivial facts about
analytic varieties in $\bC^n$ in order to compare multiplicities
of maps in two isomorphic algebras.
The conclusion is that we recover the dynamics up to piecewise
conjugacy.

For the converse, we would like to show that piecewise conjugacy
implies isomorphism of the algebras.
 
For the tensor algebra we show that the converse holds for any type of isomorphism
provided that either $n\le 3$
or the covering dimension of $X$ is at most 1.
We conjecture that this holds in complete generality.
This conjecture is backed up by the analysis of the $n=3$ case,
in which we require non-trivial topological information about
the Lie group $SU(3)$.
The conjectured converse reduces to a question about the unitary
group $U(n)$.
While the topology of $SU(n)$ and $U(n)$ gets increasing complicated for
$n\ge4$, there is reason to hope that there is a positive answer
in full generality.

On the other hand, little is
known about the converse for semicrossed products.
In Example \ref{Ex:pwc but not iso} we show that unlike the tensor algebra situation,
there are multisystems on totally a disconnected space which are piecewise conjugate
and yet their semicrossed products are not completely isometrically isomorphic. We do
not know however whether this failure can occur at the algebraic isomorphism level.

In Chapter \ref{C:semisimple}, we consider another connection
between the operator algebra and the dynamical system.
We characterize when the operator algebra, either the tensor
algebra or the crossed product, is semisimple strictly in terms
of the dynamics.
In the case of a single map, the radical of the semicrossed product
has been studied  \cite{MuRad, Pet} and finally was completely
characterized by Donsig, Katavolos and Manoussos \cite{DKM}
using a generalized notion of recurrence.
Here we introduce a notion of recurrence and
wandering sets for a dynamical system which is appropriate for
a non-commutative multivariable setting such as ours.
The main result of this chapter is the characterization of semisimplicity
in these terms.

Finally in the last chapter of this paper we mention some open problems and
further direction for future research.

\section{The one variable case}\label{S:one}

There is a long history of associating operator algebras to dynamical
systems, going back to the work of von Neumann in the 1930's.
In the self-adjoint context, one is dealing with
a (generally amenable) group of homeomorphisms.
The abstract notion of a crossed product of a C*-algebra by
an automorphism (or group of automorphisms) is an important
general construction.
There is a rich history of associating C*-invariants with the
associated dynamical systems.

The use of a nonself-adjoint operator algebras to encode a dynamical
system was first introduced by Arveson  \cite{Arvmeas}  and
Arveson--Josephson \cite{AJ} for a one-variable system $(X,\sigma)$.
In their context, $\sigma$ was a single homeomorphism with special properties.
A concrete representation was built from an appropriate invariant measure.
With certain additional hypotheses, they were able to show that
the operator algebra provided a complete invariant up to conjugacy.

The abstract version of the \textit{semicrossed product} of a dynamical
system $(X,\sigma)$ was introduced by Peters \cite{Pet}.
He does not require the existence of good invariant measures;
nor does he require $\sigma$ to be a homeomorphism.
He does require $X$ to be compact.
With this new algebra, Peters was able to show that
the semicrossed product is a complete
invariant for the system up to conjugacy
assuming that $\sigma$ has no fixed points.

In \cite{HH1}, Hadwin and Hoover considered a rather general class
of conjugacy algebras associated to a single dynamical system.
Their proofs work in considerable generality, but the semicrossed
product remains the only natural choice for the operator algebra
of a system.
Their methods allowed a further weakening of the hypotheses.
The set $\{ x \in X : \sigma^2(x)=\sigma(x) \ne x \}$ should have no interior,
but there is no condition on fixed points.
Then again they were able to recover the dynamics, up to conjugacy,
from the operator algebra.
 
In \cite{DKisotop}, we used additional information available from
studying the 2-dimensional nest representations of the semicrossed product
to completely eliminate the extraneous hypotheses on $(X,\sigma)$.
We now know that the semicrossed product, even as an algebra
without the norm structure, encodes the system up to conjugacy.
We were also able to replace a compact $X$ with a locally compact space and, as a bonus,
we were able to classify crossed products of the disc algebra by endomorphisms.

\section{Universal Operator Algebras}\label{S:universal}
 
We now discuss the choice of an appropriate covariance algebra
for the multivariable dynamical system $(X,\sigma)$.
An operator algebra encoding $(X,\sigma)$ should contain $\rC_0(X)$
as a C*-subalgebra, and there should be  $n$ elements $\fs_i$
satisfying the \textit{covariance relations}
\[ f\fs_i = \fs_i(f\circ\sigma_i) \qfor f \in \rC_0(X) \AND 1 \le i \le n.\]
This relation shows that
$
 \fs_{i_k}f_k \fs_{i_{k-1}} f_{k-1} \dots \fs_{i_1} f_1 = \fs_w g
$
where we write $\fs_w = \fs_{i_k}\fs_{i_{k-1}} \dots \fs_{i_k}$ and
$g$ is a certain product of the $f_j$'s composed with functions built from the $\sigma_i$'s.
Thus the set of polynomials in $\fs_1,\dots,\fs_n$ with coefficients in $\rC_0(X)$
forms an algebra which we call the \textit{covariance algebra} $\A_0(X,\sigma)$.
The universal algebra should be the (norm-closed non-selfadjoint)
operator algebra obtained by completing the covariance algebra in
an appropriate operator algebra norm.

Observe that in the case of compact $X$, $\A_0(X,\sigma)$ is unital,
and will contain the elements $\fs_i$ as generators.
When $X$ is not compact, it is generated by $\rC_0(X)$ and
elements of the form $\fs_i f$ for $f \in \rC_0(X)$.

By an operator algebra, we shall mean an algebra which is
completely isometrically isomorphic to a subalgebra of $\B(\H)$
for some Hilbert space $\H$.
By the Blecher--Ruan--Sinclair Theorem \cite{BRS}, there is
an abstract characterization of such algebras.
See \cite{BL, Pau} for a thorough treatment of these issues.
Our algebras are sufficiently concrete that we will not need
to call upon these  abstract results.
Nevertheless, it seems more elegant to us to define
universal operator algebras abstractly rather than in terms
of specific representations.

An operator algebra claiming to be \textit{the} operator
algebra of the system must be universal in some way.
This requires a choice of an appropriate norm condition on the generators.
A few natural choices are:
\pagebreak[3]
\begin{enumerate}
\item \textbf{Contractive:} $ \|\fs_i\| \le 1$ for $1 \le i \le n$.
\item \textbf{Isometric:} $\fs_i^* \fs_i = I$ for $1 \le i \le n$.
\item \textbf{Row Contractive:} $\big\| \begin{bmatrix}
 \fs_1 & \fs_2 & \dots & \fs_n \end{bmatrix} \big\| \le 1$.
\item \textbf{Row Isometric:} $\begin{bmatrix}
 \fs_1 & \fs_2 & \dots & \fs_n \end{bmatrix}$ is an isometry;\\
 i.e. $\fs_i^* \fs_j = \delta_{ij}$ for $1 \le i,j \le n$.
\end{enumerate}
One could add variants such as unitary, co-isometric, column contractive, etc.

In the one variable case, all of these choices are equivalent.
Indeed, the Sz.Nagy isometric dilation of a contraction is
compatible with extending the representation of $\rC_0(X)$.
This leads to the semi-crossed product introduced by Peters \cite{Pet}.
Various non-selfadjoint algebras associated to a dynamical system
(with one map) have been studied
 \cite{Arvmeas, AJ, MM, HH1, Pow, MuRad, DKM}.

Once one goes to several variables, these notions are
distinct, even in the case of commutative systems.
For example, with three or more commuting variables,
examples of Varopoulos \cite{Var} and Parrott \cite{Par}
show that three commuting contractions need not dilate
to three commuting isometries.
However a dilation theorem of Drury \cite{Dru} does
show that a strict row contraction of $n$ commuting operators
dilates to (a multiple of) Arveson's $d$-shift \cite{Arv3}.
While this is not an isometry, it is the appropriate universal
commuting row contraction.

For non-commuting variables, where there is no constraint
such as commutativity, one could dilate the $n$ contractions
to isometries separately.
We shall see that this can be done while extending
the representation of $\rC_0(X)$ to maintain the
covariance relations.
Also for the row contraction situation, there is the
dilation theorem of Frahzo--Bunce--Popescu \cite{Fra,Bun,Pop_diln}
which allows dilation of any row contraction to a row isometry.
Again we shall show that this can be done while extending
the representation of $\rC_0(X)$ to preserve the
covariance relations.

\begin{defn}
A locally compact Hausdorff space $X$ together with
$n$ proper continuous maps $\sigma_i$ of $X$ into itself
for $1 \le i \le n$ will be denoted by  $(X,\sigma)$.
We shall refer to this as a \textit{multivariable dynamical system}.
It will be called \textit{metrizable} if $X$ is metrizable.
\end{defn}

We now define the two universal operator algebras which
we will associate to $(X,\sigma)$.
We justify the nomenclature below.

\begin{defn}
Given a multivariable dynamical system $(X,\sigma)$, define the
\textbf{tensor algebra} to be the universal operator algebra
$\A(X,\sigma)$ generated by $\rC_0(X)$ and generators $\fs_1,\dots,\fs_n$
satisfying the covariance relations
\[ f\fs_i = \fs_i(f\circ\sigma_i) \qfor f \in \rC_0(X) \AND 1 \le i \le n\]
and satisfying the row contractive condition
$ \big\| \begin{bmatrix} \fs_1 & \fs_2 & \dots & \fs_n \end{bmatrix} \big\| \le 1 .$

Similarly, we define the \textbf{semicrossed product}
to be the universal operator algebra $\rC_0(X) \times_\sigma \Fn$
generated by $\rC_0(X)$ and generators $\fs_1,\dots,\fs_n$
satisfying the covariance relations and satisfying
the contractive condition $\|\fs_i\| \le 1$ for $1 \le i \le n$.
\end{defn}

We will not belabour the set theoretic issues in defining a universal
object like this, as these issues are familiar.
Suffice to say that one can fix a single Hilbert space of
sufficiently large dimension, say $\aleph_0 |X|$, on which
we consider representations of $\rC_0(X)$ and the covariance relations.
Then one puts the abstract operator algebra structure on
$\A_0(X,\sigma)$ obtained by taking the supremum over
all (row) contractive representations.
Alternatively, one forms the concrete operator algebra by taking
a direct sum over all such representations on this fixed space.

A case can be made for preferring the row contraction condition,
based on the fact that this algebra is related to other algebras
which have been extensively studied in recent years.
If $X$ is a countable discrete set, then the row contractive
condition yields the graph algebra of the underlying
directed graph that forgets which map $\sigma_i$ is
responsible for a given edge from $x$ to $\sigma_i(x)$.
In the general case, this turns out to be a C*-correspondence
algebra, or tensor algebra, as defined by Muhly and Solel \cite{MS2}.
It is for this reason that we call this algebra the tensor algebra
of the dynamical system.
As such, it sits inside a related Cuntz--Pimsner C*-algebra \cite{Pim},
appropriately defined and studied by Katsura \cite{Ka}
building on an important body of work by Muhly and Solel
beginning with \cite{MS2,MS1}.
This Cuntz-Pimsner algebra turns out to be the C*-envelope
\cite{Arv1,Ham} of the tensor algebra \cite{MS2,FMR,KK}.
The C*-envelope of the tensor algebra is therefore always nuclear.

We may consider the dynamical system $(X,\sigma)$
as an action of the free semigroup $\Fn$.
The free semigroup $\Fn$ consists of all words in the alphabet
$\{1,2,\dots,n\}$ with the empty word $\mt$ as a unit.
For each $w = i_k i_{k-1} \dots i_1$ in $\Fn$, let $\sigma_w$ denote the map
$\sigma_{i_k} \circ \sigma_{i_{k-1}}\circ \dots \circ\sigma_{i_1}$.
This semigroup of endomorphisms of $X$ induces
a family of endomorphisms of $\rC_0(X)$ by  $\alpha_w(f) = f \circ \sigma_w$.
The map taking $w \in \Fn$ to $\alpha_w$ is an antihomomorphism
of $\Fn$ into $\End(\rC_0(X))$; i.e.\ $\alpha_v \alpha_w = \alpha_{wv}$ for $v,w \in \Fn$.

This leads us to consider the contractive condition, which
is the same as considering contractive covariant representations of
the free semigroup.
Hence we call the universal algebra the semi-crossed product
$\rC_0(X) \times_\sigma \Fn$ of the dynamical system.
It also has good properties.
However we do not find this algebra as tractable
as the tensor algebra. Indeed, several problems that are
resolved in the tensor algebra case remain open
for the semicrossed product.
In particular, it often occurs (see Proposition~\ref{imposs}) that the C*-envelope
of the semicrossed product is not nuclear.

In both cases, the (row) contractive condition turns out to be equivalent
to the (row) isometric condition.  This is the result of dilation theorems
to extend (row) contractive representations to (row) isometric ones.
These are analogues of a variety of well-known dilation theorems.
The tensor algebra case is easier than the semicrossed product,
and in addition, there is a nice class of basic representations in
this case that determine the universal norm.
Indeed we exhiibit sufficiently many boundary representations to
explicitly represent the C*-envelope.
In the case of the crossed product, one needs to introduce the notion
of a \textit{full isometric dilation}; and these turn out to yield the
maximal representations of the C*-envelope.

\chapter{Dilation Theory}\label{C:diln}
\section{Dilation for the Tensor Algebra}\label{S:dilation}

We first consider a useful family of representations for the tensor
algebra analogous to those used by Peters \cite{Pet} to define the
semi-crossed product of a one variable system.

By Fock space, we mean the Hilbert space $\Fock$ with
orthonormal basis  $\{ \xi_w : w \in \Fn \}$.
This has the standard left regular representation of
the free semigroup $\Fn$ defined by
\[ L_v\xi_w = \xi_{vw} \qfor v,w \in \Fn .\]

Consider the following \textit{orbit representations} of $(X,\sigma)$.
Fix $x$ in $X$.
The orbit of $x$ is $\O(x) = \{ \sigma_w(x) : w \in \Fn \}$.
To this, we identify a natural representation of $\A(X,\sigma)$.
Define a $*$-representation $\pi_x$ of $\rC_0(X)$ on the Fock space $\F_x = \Fock$ by
$\pi_x(f) = \diag( f(\sigma_w(x)) )$, i.e.
\[  \pi_x(f) \xi_w = f(\sigma_w(x)) \xi_w \qfor f \in \rC_0(X) \AND w \in \Fn  .\]
Send the generators $\fs_i$ to $L_i$, and
let $L_x = \begin{bmatrix}L_1 & \dots & L_n \end{bmatrix}$.
Then $(\pi_x,L_x)$ is easily seen to be a covariant representation.

Define the \textit{full Fock representation} to be the (generally non-sep\-ar\-able)
representation $(\Pi,\mathbf{L})$ where $\Pi = \sum_{x \in X}^\oplus \pi_x$
and $\mathbf{L} = \sum_{x \in X}^\oplus L_x$ on $\F_X = \sum_{x \in X}^\oplus \F_x$.
We will show that the norm closed algebra generated by $\Pi(\rC_0(X))$
and $\Pi(\rC_0(X))\mathbf{L}_i$ for $1 \le i \le n$ is completely isometric
to the tensor algebra $\A(X,\sigma)$.
When $X$ is separable, a direct sum over a countable dense subset
of $X$ will yield a completely isometric copy on a separable space.

Now we turn to the dilation theorem, which is straight-forward
given our current knowledge of dilation theory.
When the dynamical system is surjective, this is closely related
to \cite[Theorem 3.3]{MS2}.

\begin{thm}\label{T:rowdiln}
Let $(X,\sigma)$ denote a multivariable dynamical system.
Let $\pi$ be a  $*$-representation of $\rC_0(X)$
on a Hilbert space $\H$, and let
$A = \begin{bmatrix} A_1 & \dots & A_n \end{bmatrix}$
be a row contraction satisfying the covariance relations
\[ \pi(f) A_i = A_i \pi(f \circ \sigma_i) \qfor 1 \le i \le n . \]
Then there is a Hilbert space $\K$ containing $\H$,
a $*$-representation $\rho$ of $\rC_0(X)$ on $\K$
and a row isometry $\begin{bmatrix} S_1& \dots & S_n \end{bmatrix}$
such that
{\renewcommand{\labelenumi}{(\roman{enumi}) }
\begin{enumerate}
\item $\rho(f) S_i = S_i \rho(f\circ\sigma_i)$ for $f \in \rC_0(X)$ and $1 \le i \le n$.
\item $\H$ reduces $\rho$ and $\rho(f)|_\H = \pi(f)$ for $f \in \rC_0(X)$.
\item $\H^\perp$ is invariant for each $S_i$, and $P_\H S_i |_\H = A_i$
for $1 \le i \le n$.
\end{enumerate}
}
\end{thm}

\begin{proof} The dilation of $A$ to a row isometry $S$ is achieved by the
Frahzo--Bunce--Popescu dilation \cite{Fra,Bun,Pop_diln}.
Consider the Hilbert space $\K = \H \otimes \Fock$ where
we identify $\H$ with $\H \otimes \bC \xi_\mt$.
Following Bunce, consider $A$ as an operator in $\B(\H^{(n)},\H)$.
Using the Schaeffer form of the isometric dilation, we can write
$D = (I_\H \otimes I_n - A^*A)^{1/2}$ in $\B(\H^{(n)})$ and
$I_\H \otimes L = \begin{bmatrix}I_\H \otimes L_1
& \dots & I_\H \otimes L_n \end{bmatrix}$.
We make the usual observation that $(\bC \xi_\mt)^\perp$
is identified with $\Fock^{(n)}$ in such a way that
$L_i|_{(\bC \xi_\mt)^\perp} \simeq L_i^{(n)}$ for $1 \le i \le n$.

Then a (generally non-minimal) dilation is obtained as
\[
 S = \begin{bmatrix} A & 0 \\JD & I_\H \otimes L^{(n)} \end{bmatrix}
 \]
 where $J$ maps $\H^{(n)}$ onto $\H \otimes \bC^n \subset \K$
 where the $i$th standard basis vector $e_i$ in $\bC^n$ is sent to $\xi_i$.
 Then
 \[ S_i = \begin{bmatrix} A_i & 0 \\JD_i & I_\H \otimes L_i^{(n)} \end{bmatrix}  \]
 where $D_i = D|{\H \otimes \bC e_i }$ is considered as an
 element of $\B(\H,\H^{(n)})$.
 
To extend $\pi$, define a $*$-representation $\rho$ on $\K$ by
\[ \rho(f) = \diag(\pi(f \circ \sigma_w) ) . \]
That is,
\[
 \rho(f) (x \otimes \xi_w) = \pi(f \circ \sigma_w)x \otimes \xi_w
 \qfor x \in \H,\ w \in \Fn .
\]
The restriction $\rho_1$ of $\rho$ to $\H \otimes  \bC^n$
is just $\rho_1(f) = \diag(\pi(f \circ \sigma_i))$.
The covariance relations for $(\pi,A)$ may be expressed as
\[ \pi(f) A = A \rho_1(f) .\]
From this it follows that $\rho_1(f)$ commutes with $A^*A$
and thus with $D$.
In particular, $\rho_1(f) D_i = D_i \pi(f\circ\sigma_i)$.
The choice of $J$ then ensures that
\[
 \rho(f) S_i |_{\H \otimes \bC \xi_\mt} =
  S_i |_{\H \otimes \bC \xi_\mt} \pi(f\circ \sigma_i) .
\]
But the definition of $\rho$ shows that
\[
  \rho(f) (I_\H \otimes L_i)
  =( I_\H \otimes L_i) \rho(f \circ \sigma_i)
\]
Hence, as $S_i$ agrees with $I_\H \otimes L_i$ on
$\H^\perp = \H \otimes (\bC \xi_\mt)^\perp$,
we obtain
\[
 \rho(f) S_i |_{\H^\perp} =
  S_i  \pi(f\circ \sigma_i)|_{\H^\perp} =
  S_i |_{\H^\perp} \pi(f\circ \sigma_i)|_{\H^\perp} .
\]
Combining these two identities yields the desired
covariance relation for $(\rho,S)$.

The other properties of the dilation are standard.
\end{proof}

\begin{rem}
If one wishes to obtain the minimal dilation, one
restricts to the smallest subspace containing $\H$
which reduces $\rho$ and each $S_i$.
The usual argument establishes uniqueness.
\end{rem}

\begin{cor}
Every row contractive representation of the covariance algebra
dilates to a row isometric representation.
\end{cor}

We now relate this to the orbit representations.
It was an observation of Bunce \cite{Bun} that the dilation $S$
of $A$ is pure if $\|A\| = r < 1$, where pure means that $S$ is
a multiple  of the left regular representation $L$.
In this case, the range $\N_0$ of the projection
$P_0 = I - \sum_{i=1}^n S_i S_i^*$ is a cyclic subspace for $S$.

Observe that for any $f \in \rC_0(X)$,
\begin{align*}
 \rho(f) S_i S_i^* &= S_i \rho(f\circ \sigma_i) S_i^*
 = S_i \big( S_i \rho(\ol{f}\circ \sigma_i)  \big)^* \\
 &= S_i  \big( \rho(f) S_i \big)^*  = S_iS_i^* \rho(f) .
\end{align*}
So $P_0$ commutes with $\rho$.
Define a $*$-representation of $\rC_0(X)$ by
$\rho_0(f) = \rho(f) |_{\N_0}$.

Then we can recover $\rho$ from $\rho_0$ and
the covariance relations.  Indeed,
$\K = \sum_{w \in \Fn}^\oplus \N_w$ where $\N_w = S_w \N_0$.
We obtain
\begin{align*}
 \rho(f) P_{\N_w} &= \rho(f) S_w P_0 S_w^*
 = S_w \rho(f\circ \sigma_w) P_0 S_w^* \\
 &= S_w \rho_0(f\circ \sigma_w) S_w^* .
\end{align*}
The spectral theorem shows that $\rho_0$ is, up to
multiplicity, a direct integral of point evaluations.
Thus it follows that the representation $(\rho,S)$
is, in a natural sense, the direct integral of
the orbit representations.
Thus its norm is dominated by the norm
of the full Fock representation.

As a consequence, we obtain:

\begin{cor} \label{C:full Fock}
The full Fock representation is a faithful completely isometric
representation of the tensor algebra $\A(X,\sigma)$.
\end{cor}

\begin{proof} By definition, if $T = \sum_{w\in\Fn} \fs_w f_w$ belongs
to $\A_0(X,\sigma)$ (i.e.\ $f_w = 0$ except finitely often),
its norm in $\A(X,\sigma)$ is determined as
\[ \|T\|_\sigma := \sup \Big\| \sum_{w\in\Fn} A_w \pi(f_w) \Big\| \]
over the set of all row contractive representations $(\pi,A)$.
Clearly, we can instead sup over the set $(\pi,rA)$ for $0 < r < 1$;
so we may assume that $\|A\| = r < 1$.
Then arguing as above, we see that $(\pi,A)$ dilates to
a row isometric representation $(\rho,S)$ which is a direct
integral of orbit representations.
Consequently the norm
\[
 \Big\| \sum_{w\in\Fn} A_w \pi(f_w)  \Big\| \le
 \Big\| \sum_{w\in\Fn} S_w \rho(f_w) \Big\| \le
 \Big\| \sum_{w\in\Fn} \mathbf{L}_w \Pi(f_w) \Big\| .
\]
Thus the full Fock representation is completely isometric,
and in particular is faithful.
\end{proof}

\begin{rem}
Indeed, the same argument shows that a faithful
completely isometric representation is obtained whenever
$\rho_0$ is a faithful representation of $\rC_0(X)$.
Conversely a representation $\rho_0$ on $\H$ induces a
Fock representation $\rho$ on $\K = \H \otimes \Fock$
by $\rho(f) = \diag(\rho_0(f \circ \sigma_w))$.
Then sending each $\fs_i$ to $I_\H \otimes L_i$ yields
a covariant representation which is faithful if $\rho_0$ is.
\end{rem}

\section{Boundary Representations and the C*-envelope}\label{S:boundary}

As mentioned in the Introduction, we are interested in the maximal dilations.
A completely contractive representation $\rho$ of an operator algebra
$\A$ on a Hilbert space $\H$ is maximal if,
whenever $\pi$ is a completely contractive dilation of $\rho$
on a Hilbert space $\K = \H \oplus \K_1$,
then $\H$ reduces $\pi$, whence $\pi$ decomposes as $\pi = \rho \oplus \pi_1$.
Such representations have the unique extension property:
if we consider $\A$ as a subalgebra of a C*-algebra $\fA$ generated
by any completely isometric image of $\A$, then
there is a unique completely positive extension of $\rho$ to $\fA$,
and it is a $*$-representation.
Such representations always factor through the C*-envelope, $\cenv(\A)$.
In particular, if one has a completely isometric maximal representation
$\rho$ of $\A$, then $\cenv(\A) = \ca(\rho(\A))$.
The maximal representations which are irreducible (no reducing subspaces)
are called boundary representations.

It follows from \cite{Arv_choq} that there are sufficiently many boundary
representations; so that their direct sum yields a completely isometric
representation of $\A$, producing the C*-envelope.
We will exhibit such representations explicitly for $\A(X,\sigma)$.

\begin{lem} \label{L:maximal}
Suppose that $\rho$ is a completely contractive representation
of $\A(X,\sigma)$ such that each $S_i = \rho(\fs_i)$ is an isometry and
\[
 \sum_{i=1}^n S_i S_i^* =
 E_\rho\big( \bigcup_{i=1}^n \sigma_i(X) \big) ,
\]
where $E_\rho$ denotes the spectral measures associated to
$\rho(\rC_0(X))$.  Then $\rho$ is maximal.
\end{lem}

\begin{proof}
Let $\rho$ be a representation of $\A(X,\sigma)$ on a Hilbert space $\H$;
and suppose that $\pi$ is any dilation on a space $\K = \H \oplus \K_1$.
The restriction of $\pi$ to $\rC_0(X)$ is a $*$-representation.
As $\H$ is invariant, it must reduce $\pi(\rC_0(X))$.
So it suffices to show that $\H$ also reduces each $\pi(\fs_i)$.

Since $f \fs_i = \fs_i (f \circ \sigma_i)$, it follows that $f \fs_i = 0$
whenever $f$ vanishes on $\sigma_i(X)$.  Therefore
$\pi(\fs_i) = E_\pi(\sigma_i(X)) \pi(\fs_i) $.
Since these isometries have orthogonal range, we always have
\[
 \sum_{i=1}^n \pi(\fs_i) \pi(\fs_i)^* \le
 E_\pi \big( \bigcup_{i=1}^n \sigma_i(X) \big) .
\]
By hypothesis,
\[
 \sum_{i=1}^n \rho(\fs_i) \rho(\fs_i)^* =
 E_\rho\big( \bigcup_{i=1}^n \sigma_i(X) \big) ,
\]
Therefore $\pi(\fs_i) \K_1$ will be orthogonal to
\[
 E_\rho\big(  \bigcup_{i=1}^n \sigma_i(X) \big) \H +
 E_\pi \big( X \setminus \bigcup_{i=1}^n \sigma_i(X) \big) \K ,
\]
which contains $\H$.
Hence $\H$ is reducing.
\end{proof}

By the results of the previous section, it suffices to find irreducible
maximal dilations of the orbit representations in order to have enough
boundary representations to determine the C*-envelope.

To do this, we need to recall the classification of atomic representations
of the Cuntz and Cuntz--Toeplitz algebras from \cite{DP1}.
An atomic representation $\pi$ of $\E_n$ on a Hilbert space $\H$
with a given orthonormal basis $\{e_n\}$ is given by $n$ isometries $S_i = \pi(\fs_i)$
with orthogonal ranges which each permute the basis up to multiplication by
scalars in the unit circle.
The irreducible atomic representations of $\E_n$ split into three types:

(1) The left regular representation $\lambda$ of $\Fn$.

(2) The infinite tail representations, which are inductive limits \mbox{of $\lambda$.}
These are obtained from an infinite sequence $\bi = i_0 i_1 i_2 \dots$
in the alphabet $\{1,\dots,n\}$.
For each $s\ge0$, let $\G_s$ denote a copy of Fock space
with basis $\{ \xi^s_w : w \in \Fn \}$.
Identify $\G_s$ with a subspace of $\G_{s+1}$ via $R_{i_s}$,
where $R_j \xi^s_w = \xi^{s+1}_{wj}$.
Set $\pi_s$ to be the representation $\lambda$ on $\G_s$.
Since $R_j$ commutes with $\lambda$, we obtain
$\pi_{s+1} R_{i_s} = R_{i_s} \pi_s$.
So we may define $\pi_\bi$ to be the inductive limit of the
representations $\pi_s$.
It is clear that the sum of the ranges of $\pi_\bi(\fs_i)$ is the whole space;
so this yields a representation of $\O_n$.
$\pi_\bi$ is irreducible if and only $\bi$ is not eventually periodic;
and two are unitarily equivalent if and only if they are shift--tail equivalent,
meaning that after deleting enough initial terms from each sequence,
they then coincide.

(3) The ring representations.
These are given by a word $u = i_1 \dots i_k$ and $\lambda \in \bT$.
Let $\C_k$ be the cyclic group with $k$ elements;
and let $\K_u$ be a Hilbert space with orthonormal basis
$\{ \xi_{s,w} : s \in \C_k ,\  w \in \F_n \setminus \F_n i_s \} $.
Define a representation $\tau_{u,\lambda}$ of $\F_n$ by
\begin{alignat*}{2}
 \tau_{u,\lambda}(\fs_i) \xi_{s+1,\mt} &= \lambda \xi_{s,\mt} &\qif &i=i_s\\
 \tau_{u,\lambda}(\fs_i) \xi_{s,w} &= \xi_{s,iw} &\qif &|w|\ge1 \OR i\ne i_s.
\end{alignat*}
This representation is irreducible if and only if $u$ is primitive
(not a power of a smaller word); and another such representation
$\tau_{v,\mu}$ is unitarily equivalent if and only if
$v$ is a cyclic permutation of $u$ and $\lambda^k=\mu^k$.

To state the theorem, we need to define some analogues
that generalize the orbit representations.

(1) The first type are the orbit representations themselves.

(2) An infinite tail representation is given by an infinite sequence
$\bi = i_0 i_1 i_2 \dots$ in the alphabet $\{1,\dots,n\}$ and a
corresponding sequence of points $x_s \in X$ for $s \ge 0$
such that $\sigma_{i_s}(x_{s+1}) = x_s$.
With the setup as in (2) above, we associate each basis vector
$\xi^s_w$ with the point $x^s_w := \sigma_w(x_s)$ in $X$.
Observe that by construction, $\sigma_{wi_s}(x_{s+1}) = \sigma_w(x_s)$;
so that these points are well defined.
Define a representation $\pi_\bi$ by defining it on the $\fs_i$
as above, and setting $\pi_\bi(f) \xi^s_w = f(x^s_w) \xi^s_w$.
This is evidently the inductive limit of the orbit representations $\pi_{x_s}$;
and thus is a completely contractive representation of $\A(X,\sigma)$.

(3) A ring representation is given by a word $u = i_1 \dots i_k$,
a scalar $\lambda \in \bT$, and a set of points $x_s \in X$ for $s \in \C_k$
satisfying $\sigma_{i_s}(x_{s+1}) = x_s$.
Again we associate a point in $X$ to each basis vector $\xi_{s,w}$
by settng $x_{s,w} := \sigma_w(x_s)$.
We define the representation $\tau_{u,\lambda}$ on the $\fs_i$ as above;
and set $\tau_{u,\lambda}(f) \xi_{s,w} = f(x_{s,w}) \xi_{s,w}$.
It is routine to verify that this is a representation of $\A(X,\sigma)$.

We can now state the result we want.

\begin{thm}\label{T:boundary}
The following are all boundary representations of
the tensor algebra $\A(X,\sigma)$.
\begin{enumerate}
\item An orbit representation $\pi_x$ for a point
 $x$ in $X \setminus \bigcup_{i=1}^n \sigma_i(X)$.
\item An infinite tail representation $\pi_\bi$ given by an infinite sequence
 $\bi = i_0 i_1 i_2 \dots$ in the alphabet $\{1,\dots,n\}$ and a
 corresponding sequence of \textbf{distinct} points $x_s \in X$ for $s \ge 0$
 such that $\sigma_{i_s}(x_{s+1}) = x_s$.
\item A ring representation $\tau_{u,\lambda}$ given by a word
 $u = i_1 \dots i_k$, a scalar $\lambda \in \bT$
 and a set of \textbf{distinct} points $x_s \in X$ for $s \in \C_k$
 satisfying $\sigma_{i_s}(x_s) = x_{s+1}$.
\end{enumerate}
\end{thm}

\begin{proof}
First let us verify that that these representations are maximal.
This is immediate from Lemma~\ref{L:maximal} because
$\sum_{i=1}^n \rho(\fs_i) \rho(\fs_i)^*$ is the identity in the
last two cases, and is $I - \xi_\mt \xi_\mt^*$ in the first case.
This is the only non-trivial case, but here the hypothesis that
$x$ is not in the range of any $\sigma_i$ means that
$\xi_\mt \xi_\mt^* \le E_{\pi_x} \big( X \setminus \bigcup_{i=1}^n \sigma_i(X) \big)$.
Indeed this is an equality, as by construction, every other basis
vector corresponds to a point in the orbit of $x$; and thus lies
in the range of $E_{\pi_x} ( \sigma_i(X) )$ for some $i$.

It remains to verify that these representations are irreducible.
The first type is irreducible because the restriction
to the algebra generated by $\fs_1,\dots,\fs_n$ is the left regular
representation, and this restriction is already irreducible.

In case (3), it follows from \cite{DP1} that the projection $P$
onto the ring space $\spn\{ \xi_{s,\mt} : s \in \C_k \}$ lies in
the \wot-closed algebra generated by $\fs_1,\dots,\fs_n$;
and a fortiori, $P$ belongs to the double commutant.
The fact that the points $x_1,\dots,x_k$ are distinct means
that the projections $\xi_{s,\mt}\xi_{s,\mt}^*$ belong to
$P \tau_u(\rC_0(X)) P$.  It follows easily now that the whole
diagonal belongs to the double commutant, and that this is
an irreducible representation.

Now consider case (2).  It may be the case that the infinite word
$\bi$ is (eventually) periodic.  But the fact that the points $x_s$
are distinct will still force irreducibility.
If $\bi$ is not eventually periodic, then the restriction to
the algebra generated by $\fs_1,\dots,\fs_n$ is already irreducible.

So suppose that $\bi$ is periodic from some point on,
say repeating some word primitive $u$ of length $k$.
Then we can truncate it, and then continue it in both directions
as a periodic sequence.
It follows from \cite{DP1} that the $P$ onto the spine lies
in the double commutant.
Now $P\H$ looks like $\ltwo(\bZ)$, and $P \pi_\bi(\rC_0(X)) |_{P\H}$
is a subalgebra of the diagonal algebra which separates points in the
tail because the points $x_s$ are distinct, say for $s \ge 0$.
Moreover $P \sum_{i=1}^n \pi_\bi(\fs_i) P$ is the bilateral shift
on  $P\H$.  This does not commute with any non-scalar diagonal operator.
So the restriction to $P\H$ is \wot-dense.
It follows easily now that $\pi_\bi$ is irreducible.
\end{proof}

\begin{cor}\label{C:Cenv tensor}
The C*-envelope of $\A(X,\sigma)$ is $\ca(\pi(\A(X,\sigma))$,
where $\pi$ is the direct sum of the boundary representations
described in Theorem~$\ref{T:boundary}$.
\end{cor}

\begin{proof}
By Corollary~\ref{C:full Fock}, the full Fock representation $\Pi$ is
a completely isometric representation of $\A(X,\sigma)$.
Therefore any maximal dilation $\pi$ of $\Pi$ will yield the C*-envelope
as desired.  Thus it suffices to show that each orbit representation
$\pi_x$ dilates to one of the representations described in
Theorem~\ref{T:boundary} or to a direct sum or integral
of such representations.

Given $x \in X$, we can recursively select points $x_s \in X$ and
integers $1 \le i_s \le n$ for $s \ge 0$ such that $x_0=x$ and
$\sigma_{i_s} x_{s+1} = x_s$.
If at some point this procedure stops, it is because some $x_{s_0}$
lies in $X \setminus \bigcup_{i=1}^n \sigma_i(X)$.
In this case, it is evident that $\pi_{x_{i_{s_0}}}$ is a dilation
of $\pi_x$ because $x$ lie in the orbit of $x_{s_0}$.
This is maximal of the first type.

If the procedure can be repeated ad infinitum, then either
the points $x_s$ are all distinct or there is a repetition.
In the first case, one obtains a maximal dilation of $\pi_x$ to
$\pi_\bi$ for the sequence just constructed.

Finally, we suppose that there is a repetition $x_{s_0}=x_{s_1}$
for some $0 \le s_0 < s_1$, with $x_s$ all distinct in this range.
Let $u= i_{s_0+1}\dots i_{s_1}$.
If $x$ is not in the set $\{x_s : s_0 \le s < s_1 \}$, then
it lies in the orbit but not in the ring.
Therefore the representation $\tau_{u,\lambda}$ is a dilation
of $\pi_x$, and is maximal of the third type.

However, it may happen that $x=x_s$ for some point in the
ring and does not occur elsewhere in the orbit.
In this case, no single $\tau_{u,\lambda}$ is a dilation of $\pi_x$.
However the direct integral $\pi$ of these representations with respect
to Lebesgue measure on the circle effectively replaces the
one dimensional space associated to each vertex with
a copy of $\ltwo(\bZ)$, and replaces the
rank one operators taking each vector in the ring to the next
with copies of the bilateral shift.
Then taking any basis vector $\xi$ in the Hilbert space over $\xi_{s,\mt}$,
it is clear that the restriction of $\pi$ to the
cyclic subspace $\pi(\A(X,\sigma)) \xi$ is equivalent to $\pi_x$.
Hence $\pi$ dilates $\pi_x$.
As a direct integral of maximal dilations, it is clearly maximal.
(One can invoke Lemma~\ref{L:maximal} instead, if one wishes.)
\end{proof}

Now we can establish the converse of Lemma~\ref{L:maximal},
yielding a significant strengthening of Theorem~\ref{T:rowdiln}..

\begin{cor} \label{C:maximal}
A completely contractive representation $\rho$ of
the tensor algebra $\A(X,\sigma)$
is maximal if and only if $S_i = \rho(\fs_i)$ is an isometry
for $1 \le i \le n$ and
\[
 \sum_{i=1}^n S_i S_i^* =
 E_\rho\big( \bigcup_{i=1}^n \sigma_i(X) \big) ,
\]
where $E_\rho$ denotes the spectral measures associated to
$\rho(\rC_0(X))$.
\end{cor}

\begin{proof}
One direction follows from Lemma~\ref{L:maximal}.

Conversely, one can easily check that in each of the boundary representations
$\pi$ of Theorem~\ref{T:boundary}, each $S_i = \pi(\fs_i)$ is an isometry and
$\sum_{i=1}^n S_i S_i^* =
 E_ \pi \big( \bigcup_{i=1}^n \sigma_i(X) \big)$.
As these facts are preserved by $*$-homomorphisms,
they also hold in the C*-envelope.
Any maximal dilation extends uniquely to a $*$-representation
of $\cenv(\A(X,\sigma))$, and therefore also inherits these properties.
\end{proof}

Of course, it should be possible to construct the maximal dilations
more directly.  We explain how to do this when $X$ is metrizable.

Clearly, if $\rho$ is maximal, then each $S_i = \rho(\fs_i)$ is isometric by
Theorem~\ref{T:rowdiln}.
Observe that
\[
 \rho(f) S_i S_i^* =
 S_i \rho(f \circ \sigma_i) S_i^* =
 S_i S_i^* \rho(f) .
\]
Hence $\sum_{i=1}^n S_i S_i^*$ commutes with $\rho(\rC_0(X))$.
Suppose that there is an integer $i_0$ so that
\[
 P = E_\rho(\sigma_{i_0}(X))
       \Big(  \sum_{i=1}^n \rho(\fs_i) \rho(\fs_i)^* \Big)^\perp \ne 0 .
\]

Since $X$ is metrizable, we can
choose a Borel selector $\omega$ for $\sigma_{i_0}^{-1}$;
that is, a Borel function $\omega$ taking $\sigma_{i_0}(X)$ to $X$
so that $\sigma_i \circ \omega_i = \id_{\sigma_i(X)}$.
The existence of such a function is elementary \cite[Theorem~4.2]{Part}.

Let $\H_0$ be the Hilbert space $P\H$;
and let $J$ denote the natural injection of $\H_0$ into $\H$.
Define a representation $\pi_0$ of $\rC_0(X)$ on $\H_0$ by
$\pi_0(f) = P \rho(f \circ \omega)|_{P\H}$.
Now define a dilation of $\rho$ on the space $\K = \H_0 \oplus \H$ by
$\pi(f) = \pi_0(f) \oplus \rho(f)$ and
\[
  \pi(\fs_{i_0}) = \begin{bmatrix}0&0\\ J&S_{i_0} \end{bmatrix}
  \qand
  \pi(\fs_i) = \begin{bmatrix}0&0\\0&S_i \end{bmatrix} \FOR i \ne i_0.
\]
It is easy to verify that this is a bona fide dilation which clearly
does not decompose as a direct sum.  Hence $\rho$ is not maximal.

This dilation may be dilated to an isometric one by Theorem~\ref{T:rowdiln};
and this pair of operations may be repeated until a maximal dilation is obtained.

\section{C*-correspondences}\label{S:correspondence}

In this section, we identify the tensor algebra with the tensor algebra
of a C*-correspondence in the sense of Pimsner \cite{Pim},
Muhly--Solel \cite{MS2} and Katsura \cite{Ka}.
This will provide another description of the C*-envelope of $\A(X,\sigma)$.

Define $E = X \times n$ to be the union of $n$ disjoint copies of the space $X$.
We will view $\E = \rC_0(E)$ as a C*-correspondence over $\rC_0(X)$.
To this end, we need to define left and right actions of $\rC_0(X)$
and define a $\rC_0(X)$-valued inner product.
Let $\xi = \xi(x,j)$ and $\eta$ denote elements of $\E$,
and $f \in \rC_0(X)$.
The actions are given by
\begin{align*}
(\xi \cdot f)(x,j) &= \xi(x,j) f(x)  \\
(f  \cdot \xi)(x,j) &:= \phi(f) \xi(x,j) = \xi(x,j) (f\circ\sigma_j) (x) \\
\intertext{and the inner product is}
\lip \xi \mid \eta \rip (x) &= \sum_{j=1}^n \ol{\xi(x,j)} \eta(x,j) .
\end{align*}
That is, $\E$ is a (right) Hilbert C*-module over $\rC_0(X)$  \cite{L}
with the additional structure as a left module over $\rC_0(X)$.

As a Hilbert C*-module, $\E$ has the operator space
structure of column $n$-space $\operatorname{Col}_n(\rC_0(X))$.
The adjointable left multipliers of $\E$
form a C*-algebra $\fL(\E)$.
It will be convenient to explicitly identify the left action as
a $*$-homomorphism $\phi$ of $\rC_0(X)$ into $\fL(\E)$ by
\[ \phi(f) = \diag(f\circ\sigma_1,\dots,f\circ\sigma_n) .\]
We will usually write $\phi(f) \xi$ rather than $f\xi$ to emphasize the
role of $\phi$, as is common practice.

In general, $\phi$ is not faithful.
Indeed, set $U_0 = X \bsl \bigcup_{i=1}^n \sigma_i(X)$
which is open because the $\sigma_i$ are proper.
Then $\ker\phi$ consists of all functions with support
contained in $U_0$.
However the Hilbert C*-module is full, since $\lip \E,\E \rip = \rC_0(X)$.
Also the left action is essential, i.e.\ $\phi(\rC_0(X) ) \E = \E$.

We briefly review Muhly and Solel's construction of the
tensor algebra of $\E$ and two related C*-algebras.
Set $\E^{\otimes 0} = \rC_0(X)$ and
\[ \E^{\otimes k} = \underbrace{
\E \otimes_{\rC_0(X)} \E \otimes_{\rC_0(X)} \dots \otimes_{\rC_0(X)} \E
}_{k \text{ copies}} \qfor k\ge1 .
\]
Notice that $\xi f \otimes \eta = \xi \otimes \phi(f) \eta$.

Let $\ep_i$ denote the column vector with a $1$ in the $i$th position.
A typical element of $\E$ has the form $\sum_{i=1}^n \ep_i f_i$
for $f_i \in \rC_0(X)$.
For each word $w = i_k \dots i_1 \in \Fn$, write
\[ \ep_w := \ep_{i_k} \otimes \dots \otimes \ep_{i_1} .\]
A typical element of $\E^{\otimes k}$ has the form $\sum_{|w|=k} \ep_w f_w$
for $f_w \in \rC_0(X)$.

Naturally, $\E^{\otimes k}$ is a $\rC_0(X)$-bimodule with the rules
\[ (\xi_k \otimes \dots \otimes \xi_1) \cdot f = \xi_k \otimes \dots \otimes (\xi_1 f) \]
and
\[
f \cdot  (\xi_k \otimes \dots \otimes \xi_1) = (\phi(f) \xi_k) \otimes \dots \otimes \xi_1 .
\]
Observe that
\begin{align*}
 f \cdot \ep_w &= (\phi(f) \ep_{i_k}) \otimes  \ep_{i_{k-1}} \otimes \dots \otimes \ep_{i_1}\\
 &= \ep_{i_k} (f\circ \sigma_{i_k}) \otimes  \ep_{i_{k-1}} \otimes\dots \otimes \ep_{i_1} \\
 &= \ep_{i_k} \otimes \ep_{i_{k-1}} (f\circ \sigma_{i_k}\circ\sigma_{i_{k-1}}) \otimes\dots \otimes \ep_{i_1} \\
 &= \ep_{i_k} \otimes  \ep_{i_{k-1}} \otimes\dots \otimes
  \ep_{i_1}( f\circ \sigma_{i_k}\circ\sigma_{i_{k-1}}\circ\cdots\circ\sigma_{i_1}) \\
 &=  \ep_w  (f\circ \sigma_w) .
\end{align*}
This identifies a $*$-homomorphism $\phi_k$ of $\rC_0(X)$ into
$\fL(\E^{\otimes k})$ by
\[  \phi_k(f) = \diag(f \circ \sigma_w)_{|w|=k}  , \]
namely
\[ \phi_k(f) \sum_{|w|=k}  \ep_w g_w = \sum_{|w|=k}  \ep_w (f \circ\sigma_w) g_w .\]

The inner product structure is defined recursively by the rule
\[
 \ip{ \xi \otimes \mu, \eta \otimes \nu} =
 \ip{ \mu, \phi(\ip{\xi ,\eta }) \nu } \qforal \xi,\eta \in \E,\ \mu,\nu \in \E^{\otimes k} .
\]
This seems complicated, but in our basis it is transparent:
\[
 \big\lip \sum_{|w|=k} \ep_w f_w,  \sum_{|w|=k} \ep_w g_w \big\rip
 =  \sum_{|w|=k} \ol{f}_w g_w .
\]

The Fock space of $\E$ is $\F(\E) = \sum_{n\ge0}^\oplus \E^{\otimes n}$.
This becomes a C*-correspondence as well,
with the  $\rC_0(X)$-bimodule actions already defined on each summand,
and the $\rC_0(X)$-valued inner product obtained by declaring the
summands to be orthogonal.
In particular, this yields a $*$-isomorphism $\phi_\infty$
of $\rC_0(X)$ into $\fL(\F(\E))$ by $\phi_\infty(f)|_{ \E^{\otimes k}}= \phi_k(f)$.

There is a natural tensor action of $\xi \in \E$
taking $\E^{\otimes k}$ into $\E^{\otimes k+1}$ given by
\[
 T_\xi^{(k)} (\xi_1 \otimes \dots \otimes \xi_k) = \xi \otimes \xi_1 \otimes \dots \otimes \xi_k .
\]
Define $T_\xi$ acting on $\F(\E)$ by setting $T_\xi|_{\E^{\otimes k}} = T_\xi^{(k)}$.
The tensor algebra $\T_+(\E)$ of $\E$ is the norm-closed non-selfadjoint
subalgebra generated by $\phi_\infty(\rC_0(X))$ and $\{T_\xi : \xi \in \E \}$.
The C*-algebra generated by $\T_+(\E)$ is called the Toeplitz C*-algebra $\T(\E)$.

We will show that $\T_+(\E)$ is completely isometrically isomorphic to $\A(X,\sigma)$.
In addition, we shall show that the enveloping C*-algebra of $\A(X,\sigma)$
in the full Fock representation is $*$-isomorphic to $\T(\E)$.
Moreover, we will identify the quotient of $\T(\E)$ which is the
C*-envelope of $\A(X,\sigma)$.
This will be an application of Katsura \cite{Ka}  and Katsoulis and Kribs \cite{KK3}.
To describe this quotient, we need some further description of the
work of Katsura.

As usual, the set $\fK(\E)$ of compact multipliers is the
ideal of elements of $\fL(\E)$ generated by the rank one
elements $\theta_{\xi,\eta}$ for $\xi,\eta \in \E$ given by
$\theta_{\xi,\eta} \zeta = \xi \ip{ \eta,\zeta}$.
In our case, $\phi(\rC_0(X))$ is contained in $\fK(\E)$.
To see this, given $f\in \rC_0(X)$, factor $f \circ\sigma_i = g_i h_i$.
Then for $\zeta = \sum_{j=1}^n \ep_j k_j \in \E$,
\begin{align*}
 \sum_{i=1}^n \theta_{\ep_i g_i, \ep_i \bar{h}_i} \zeta &=
 \sum_{i=1}^n \ep_i g_i \ip{  \ep_i \bar{h}_i, \zeta} =
 \sum_{i=1}^n \ep_i g_i h_i k_i \\  &=
\sum_{i=1}^n \ep_i (f \circ\sigma_i) k_i = \phi(f) \zeta
\end{align*}

Katsura's ideal for a C*-correspondence $\E$ over a C*-algebra $\fA$ is defined as
\[ \J_\E = \phi^{-1}(\fK(\E)) \cap \operatorname{ann}(\ker \phi) \]
where $\operatorname{ann}(\I) = \{ a \in \fA : ab = 0 \FORAL b \in \I \}$.
In our setting, since $\ker \phi$ consists of functions supported on $U_0$,
$\J_\E = I(\ol{U_0})$, the ideal of functions vanishing on $\ol{U_0}$.

The space $\F(\E) \J_\E$ becomes a Hilbert C*-module over $\J_\E$.
Moreover $\fK(\F(\E)\J_\E)$ is spanned by terms of the form
$\theta_{\xi f, \eta}$ where $\xi,\eta \in \F(\E)$ and $f \in \J_\E$.
He shows that this is an ideal in $\T(\E)$.
The quotient $\O(\E) = \T(\E)/\fK(\F(\E)\J_\E)$ is the
Cuntz--Pimsner algebra of $\E$.
(Note that in Muhly--Solel \cite{MS2}, this is called the relative Cuntz-Pimsner
algebra $\O(\J_\E,\E)$.  However the crucial role of this particular ideal $\J_\E$
is due to Katsura.)

A representation of a C*-correspondence $\E$ consists of a linear map
$t$ of $\E$ into $\B(\H)$ and a $*$-representation $\pi$ of $\rC_0(X)$ on $\H$
such that
\begin{enumerate}
\item $t(\xi)^* t(\eta) = \pi(\ip{\xi,\eta}) \qforal \xi,\eta \in \E$
\item $\pi(f) t(\xi) = t(\phi(f) \xi) \qforal f \in \rC_0(X),\ \xi \in \E$.
\end{enumerate}
Such a representation is automatically a right module map as well.
Moreover, when $\pi$ is injective, $t$ is automatically an isometry.
Denote by $\ca(\pi,t)$ the C*-algebra generated by $\pi(\rC_0(X))$
and $t(\E)$.

There is a universal C*-algebra $\T_\E$ generated by such representations, and
Katsura \cite[Prop.6.5]{Ka} shows that the universal C*-algebra
is isomorphic to $\T(\E)$.

There is an induced $*$-representation of $\fK(\E)$ given by
\[ \psi_t(\theta_{\xi,\eta}) = t(\xi) t(\eta)^* .\]
Katsura shows that if $\pi(f)$ belongs to $\psi_t(\fK(\E))$,
then $f \in \J_\E$ and $\pi(f) = \psi_t(\phi(f))$.
Then he introduces an additional property of a representation
which he calls covariance.  Since we already have a different
and natural use for this term, we shall call such a
representation \textit{reduced} if it satisfies
\begin{enumerate}\addtocounter{enumi}{2}
\item $\pi(f) = \psi_t(\phi(f)) \qforal f \in \J_\E$.
\end{enumerate}
There is again a universal C*-algebra $\O_\E$
for reduced representations of $\E$.
This algebra is shown to be $*$-isomorphic to $\O(\E)$ \cite[Prop.6.5]{Ka}.

One says that $\ca(\pi,t)$ admits a gauge action if
there is a map $\beta$ of the circle $\bT$ into $\Aut(\ca(\pi,t))$
such that $\beta_z(\pi(f)) = \pi(f)$ and $\beta_z(t(\xi)) = z t(\xi)$
for $f\in\rC_0(X)$ and $\xi \in \E$.
The universal algebras have this property automatically.

Every operator algebra is contained completely isometrically in a
canonical minimal C*-algebra known as its C*-envelope.
Muhly and Solel \cite[Theorem~6.4]{MS2} show that, when $\phi$ is
injective, that the C*-envelope of $\T_+(\E)$ is $\O(\E)$.
This was further analyzed in \cite{FMR}.
In the non-injective case, this is done by Katsoulis-Kribs \cite{KK}.

\begin{thm}\label{correspondence}
Let $(X,\sigma)$ be a multivariable dynamical system, and let
$\E$ be the associated C*-correspondence.
Then $\A(X,\sigma)$ is completely isometrically isomorphic
to the tensor algebra $\T_+(\E)$.
Consequently, the C*-envelope of $\A(X,\sigma)$ is
$*$-isomorphic to $\O(\E)$.
\end{thm}

\begin{proof} The point is to observe that the Fock representations give rise
to representations of $\E$ which are sufficient to yield a faithful
representation of $\T_\E$.
Fix $x \in \X$.
Define a representation of $\E$ on $\F_x = \Fock$
by
\[
 \pi_x(f) = \diag(f(\sigma_w(x))) \qand t_x(\xi) = \sum_{i=1}^n L_i \pi_x(g_i)
\]
 for $f \in \rC_0(X)$ and $\xi = \sum_{i=1}^n \ep_i g_i \in \E$.
It is routine to verify for $\eta = \sum_{i=1}^n \ep_i  h_i$ that
\[ t_x(\xi)^* t_x(\eta) = \sum_{i=1}^n \pi_x(\ol{g}_i h_i) = \pi_x(\ip{\xi,\eta}) \]
and
\begin{align*}
\pi_x(f) t_x(\xi) &= \sum_{i=1}^n \pi_x(f) L_i \pi_x(g_i) \\
&= \sum_{i=1}^n L_i \pi_x(f\circ\sigma_i) \pi_x(g_i) = t_x(\phi(f)\xi) .
\end{align*}

Moreover it is also clear that the C*-algebra generated is
exactly the Fock representation of $\A(X,\sigma)$ for the point $x$.
Again we take the direct sum over all $x \in X$ (or a countable
dense subset in the separable case) to obtain a full Fock
representation $(\Pi,T)$.
The resulting C*-algebra admits a gauge action by conjugating
on each Fock space $\F_x$ by the unitary operator $U_z = \diag(z^{|w|})$.
Since the representation $\Pi$ is obviously faithful,
Katsura's Theorem~6.2 of \cite{Ka} shows that
we obtain a faithful representation of $\T_\E$ provided that
\[ \I_{\Pi,T} := \{ f \in \rC_0(X) : \Pi(f) \in \psi_T(\fK(\E)) \} = 0 .\]

To see this, observe that if $\xi = \sum_{i=1}^n \ep_i g_i \in \E$
and $\eta = \sum_{i=1}^n \ep_i  h_i$, then
\begin{align*}
 \psi_T(\theta_{\xi,\eta}) &= \sum_{x\in X}^\oplus t_x(\xi) t_x(\eta)^*
 = \sum_{x\in X}^\oplus \sum_{i,j=1}^n L_i g_i \ol{h}_j L_j^* .
\end{align*}
All of the vectors $\xi_{\mt,x} \in \F_x$ lie in the kernel of all
of these maps.  On the other hand, if $f \ne 0$, then $f(x) \ne 0$
for some $x$ and $\pi_x(f) \xi_{\mt,x} = f(x) \xi_{\mt,x}  \ne 0$.
This establishes the claim.

Now the rest follows from the discussion preceding the theorem.
\end{proof}

\begin{eg}
We conclude this section with an example to illuminate these ideas.
Consider $X = [0,1]$ and let $\sigma_1(x) = x/3$ and $\sigma_2(x) = (2+x)/3$.
This is an iterated function system.
Therefore there is a unique non-empty compact subset $Y$
such that $Y = \sigma_1(Y) \cup \sigma_2(Y)$.
In this case, it is easily seen to be the Cantor set $X_\infty$.

We form the algebra $\A(X,\sigma)$ and the C*-correspondence $\E$.
Observe that the kernel $\ker \phi = I\big([0,\frac13]\cup [\frac23,1]\big)$.
Thus $\J_\E = I([\frac13,\frac23])$.

Katsura shows that $\phi_\infty(\rC_0(X))$ has $\{0\}$ intersection
with the ideal $\fK(\F(\E)\J_\E)$.  Hence the representation of
$\rC_0(X)$ into $\T(\E)/\fK(\F(\E)\J_\E)$ is injective.
However $\phi_\infty(\rC_0(X))$ does intersect $\fK(\F(\E))$.
To see this, note that $\phi_\infty(f) = \diag(\phi^{(k)}(f))$ belongs to
$\fK(\F(\E))$ if and only if
\[ \lim_{k\to\infty} \| \phi^{(k)}(f) \| = 0 .\]
This is because each $\phi^{(k)}(f)$ belongs to $\fK(\F(\E))$
for every $f$.  (This was noted above for $\E$, but works just
as well for $\E^{\otimes k}$.)
Let $X_k = \bigcup_{|w|=k} \sigma_w(X)$.
This is a decreasing sequence of compact sets with
$\bigcap_{k\ge1} X_k = X_\infty$.
Thus it is easy to see that $\phi_\infty^{-1}(\fK(\F(\E))) = I(X_\infty)$.
In particular, the quotient of $\T(\E)$ by $\fK(\F(\E))$ is not injective
on $\rC_0(X)$.
\end{eg}

\section{Dilation and the Semi-crossed Product}

In this section, we consider the contractive case.
Again a routine modification of the classical theory
yields a dilation to the isometric case.
This will allow us to determine certain faithful
representations of $\rC_0(X) \times_\sigma \Fn$.

\begin{prop}\label{P:dilncont}
Let $(X,\sigma)$ denote a multivariable dynamical system.
Let $\pi$ be a  $*$-representation of $\rC_0(X)$
on a Hilbert space $\H$, and let $A_1, \dots, A_n$
be contractions satisfying the covariance relations
\[ \pi(f) A_i = A_i \pi(f \circ \sigma_i) \qfor 1 \le i \le n . \]
Then there is a Hilbert space $\K$ containing $\H$,
a $*$-representation $\rho$ of $\rC_0(X)$ on $\K$
and isometries $S_1,\dots, S_n$ such that
{\renewcommand{\labelenumi}{(\roman{enumi}) }
\begin{enumerate}
\item $\rho(f) S_i = S_i \rho(f\circ\sigma_i)$ for $f \in \rC_0(X)$ and $1 \le i \le n$.
\item $\H$ reduces $\rho$ and $\rho(f)|_\H = \pi(f)$ for $f \in \rC_0(X)$.
\item $\H^\perp$ is invariant for each $S_i$, and $P_\H S_i |_\H = A_i$
for $1 \le i \le n$.
\end{enumerate}
}
\end{prop}

\begin{proof}
This time, we can dilate each isometry separately in the classical way
so long as they use pairwise orthogonal subspaces for these extensions.
Form $\K = \H \otimes \Fock$, and again identify
$\H$ with $\H \otimes \bC \xi_\mt$.
Let $D_i = (I - A_i^*A_i)^{1/2}$.
As before, we make the identification of $(\bC \xi_\mt)^\perp$
with $\Fock^{(n)}$ in such a way that
$L_i|_{(\bC \xi_\mt)^\perp} \simeq L_i^{(n)}$ for $1 \le i \le n$.
Let $J_i = I_\H \otimes L_i|_{\bC \xi_\mt}$ be the isometry
of $\H \otimes \bC \xi_\mt$  onto $\H \otimes \bC \xi_i$ for $1 \le i \le n$.

Define a $*$-representation $\rho$ of $\rC_0(X)$ as before by
\[ \rho(f) = \diag(\pi(f \circ \sigma_w)) ; \]
and define isometric dilations of the $A_i$ by
\[
 V_i = \begin{bmatrix}
             A_i & 0 \\
            J_iD_i & I_\H \otimes L_i^{(n)}
           \end{bmatrix}
\]
Again we have identified $(\bC \xi_\mt)^\perp$ with $\Fock^{(n)}$.

To verify the covariance relations, compute:
for $x\in\H$, $f \in \rC_0(X)$ and $w \in \Fn\bsl\{\mt\}$
\begin{align*}
 \rho(f) V_i (x \otimes \xi_w)
 &= \rho(f) (x \otimes \xi_{iw}) \\
 &= \pi(f\circ\sigma_i\circ\sigma_w) x \otimes \xi_{iw} \\
 &= V_i (\pi(f\circ\sigma_i\circ\sigma_w) x \otimes \xi_w) \\
 &= V_i \rho(f\circ\sigma_i)  (x \otimes \xi_w ) .
\end{align*}
While if $w = \mt$,
\begin{align*}
 \rho(f) V_i ( x \otimes \xi_\mt )
 &= \rho(f) (A_i x \otimes \xi_\mt + D_i x \otimes \xi_i )\\
 &= (\pi(f) A_i x) \otimes \xi_\mt + (\pi(f\circ\sigma_i) D_i x) \otimes \xi_i
\end{align*}
and
\begin{align*}
 V_i \rho(f\circ\sigma_i)(x \otimes \xi_\mt)
 &= (A_i \pi(f\circ\sigma_i) x) \otimes \xi_\mt +  (D_i \pi(f\circ\sigma_i) x )\otimes \xi_i \\
 &= (\pi(f) A_i x)  \otimes \xi_\mt +  (D_i \pi(f\circ\sigma_i) x )\otimes \xi_i
\end{align*}
Thus we will have the desired relation provided that $D_i$ commutes
with $\pi(f\circ\sigma_i)$.
This is true and follows from
\[
 A_i^*A_i \pi(f\circ\sigma_i) = A_i^* \pi(f) A_i = \pi(f\circ\sigma_i) A_i^*A_i . \qedhere
\]
\end{proof}

\begin{cor}
Let $(X,\sigma)$ denote a multivariable dynamical system.
Every contractive covariant representation of $(X,\sigma)$
dilates to an isometric representation.
\end{cor}

Let $(\pi,S)$ be an isometric representation of $(X,\sigma)$.
The covariance relations extend to
the abelian von Neumann algebra $\pi(\rC_0(X))''$.
This algebra has a spectral measure $E_\pi$ defined
on all Borel  subsets of $X$.
Indeed, there is a $*$-representation $\ol{\pi}$ of
the C*-algebra $\operatorname{Bor}(X)$ of all
bounded Borel functions on $X$ extending $\pi$.
For any Borel set $A\subset X$, the covariance relations say that
\[  E_\pi(A) S_i = S_i E_\pi(\sigma_i^{-1}(A))  .\]
or equivalently
\[  E_\pi \big( \sigma_i^{-1}(A) \big) =  S_i^* E_\pi(A) S_i .\]
However,
\[ S_i E_\pi \big( \sigma_i^{-1}(A) \big) S_i^* = E_\pi(A) S_i S_i^* \le E_\pi(A) . \]
This is not completely satisfactory.

\begin{defn}
An isometric covariant representation $(\pi,S)$ is
a \textit{full isometric representation}
provided that $S_iS_i^* = E_\pi(\sigma_i(X))$.
\end{defn}

The calculation above shows that a full isometric representation
has the following important property.

\begin{lem}\label{L:full}
If $(\pi,S)$ is a full isometric covariant representation, then
\[ S_i E_\pi \big( \sigma_i^{-1}(A) \big) S_i^* = E_\pi(A) \]
for all Borel subsets $A \subset \sigma_i(X)$.
\end{lem}

The following result is a dilation in the spirit of unitary dilations,
rather than isometric dilations, in that the new space contains
$\H$ as a semi-invariant subspace (rather than a coinvariant
subspace).

Our proof requires the existence of a Borel cross section
for the inverse of a continuous map, and that appears
to require a separability condition.
We were not able to find any counterexamples for
huge spaces, but the case of most interest is, in any case,
the metrizable one.
Recall that a locally compact Hausdorff space is
metrizable if and only if it is second countable.
In this case, it is also separable.

\begin{thm}\label{T:fulldiln}
Let $(X,\sigma)$ be a metrizable multivariable dynamical system.
Every isometric covariant representation $(\pi,S)$ of $(X,\sigma)$
has a dilation to a full isometric representation in the sense that there
is a Hilbert space $\K$ containing $\H$, a $*$-representation
$\rho$ of $\rC_0(X)$ and isometries $T_1,\dots.T_n$
such that
\begin{enumerate}
\item $(\rho,T)$ is a full isometric representation of $(X,\sigma)$.
\item $\H$ reduces $\rho$ and $\rho|_\H = \pi$, and
\item $\H$ is semi-invariant for each $T_i$ and $P_\H T_i|_\H = S_i$.
\end{enumerate}
\end{thm}

\begin{proof}  Choose a Borel selector $\omega_i$ for each $\sigma_i^{-1}$;
that is, a Borel function $\omega_i$ taking $\sigma_i(X)$ to $X$
so that $\sigma_i \circ \omega_i = \id_{\sigma_i(X)}$.
The existence of such a function is elementary \cite[Theorem~4.2]{Part},
but requires metrizability.

Observe first that $S_iS_i^*$ commutes with $\pi$.
For any $f \in \rC_0(X)$,
\begin{align*}
 \pi(f) S_i S_i^* &= S_i \pi(f\circ \sigma_i) S_i^*
 = S_i \big( S_i \pi(\ol{f}\circ \sigma_i)  \big)^* \\
 &= S_i  \big( \pi(f) S_i \big)^*  = S_iS_i^* \pi(f) .
\end{align*}

The construction is recursive.
Let $\H_i = E_{\pi}(\sigma_i(X)) (I-S_iS_i^*)\H$, and let
$J_i$ be the natural injection of $\H_i$ into $\H$.
Since $\H_i$ reduces $\ol{\pi}$, we may define a
$*$-representation $\ol{\pi}_i$ as its restriction.
Form a Hilbert space $\L_1 = \H \oplus \sum_{1 \le i \le n}^\oplus \H_i$.
Define a $*$-representation $\pi_1$ on $\L_1$ by
\[ \pi_1(f) = \pi(f) \oplus \sum_{1 \le i \le n}^\oplus \ol{\pi}_i(f \circ \omega_i) \]
Also define an extension $A^{(1)}_i$ of $S_i$ by $A^{(1)}_i|_\H = S_i$,
$A^{(1)}_i|_{\H_j} = 0$ for $j \ne i$ and $A^{(1)}_i|_{\H_i} = J_i$.
Then since $f \circ \sigma_i \circ \omega_i|_{\sigma_i(X)} = f|_{\sigma_i(X)}$,
\begin{align*}
 \pi_1(f) A^{(1)}_i &= \pi_1(f) (S_i+J_i) = \pi(f) S_i + \pi(f) J_i P_{\H_i} \\
 &=S_i  \pi(f\circ\sigma_i) + J_i P_{\H_i} E_\pi(\sigma_i(X)) \pi(f) \\
 &=  S_i  \pi(f\circ\sigma_i) + J_i \ol{\pi}_i(f \circ \sigma_i \circ \omega_i) \\
 &= (S_i  + J_i) \pi_1(f \circ \sigma_i)
 = A^{(1)}_i \pi_1(f \circ \sigma_i) .
\end{align*}
So this is a contractive covariant representation of $(X,\sigma)$.
Hence it has an isometric dilation $(\rho_1,T^{(1)})$
on a Hilbert space $\K_1$.
Moreover, by construction, the range of $T^{(1)}_i$
contains $E_\pi(\sigma_i(X)) \H$.

Repeat this procedure with the representation $(\rho_1,T^{(1)})$
to obtain an isometric representation $(\rho_2,T^{(2)})$ such that
the range of $T^{(2)}_i$ contains $E_{\rho_1}(\sigma_i(X)) \K_1$.
By induction, one obtains a sequence of dilations $(\rho_k,T^{(k)})$
on an increasing sequence of Hilbert spaces $\K_k$ so that, at every stage,
the range of $T^{(k+1)}_i$ contains $E_{\rho_k}(\sigma_i(X)) \K_k$.

Let $\K$ be the Hilbert space completion of the union of the $\K_i$.
Define a $*$-representation $\rho$ of $\rC_0(X)$ on $\K$
by $\rho|_{\K_k} = \rho_k$, and set $T_i|_{\K_k} = T^{(k)}_i$.
This is a full isometric dilation by construction.
\end{proof}

We now show that the full isometric representations are
precisely the maximal representations.

Douglas and Paulsen \cite{DouPau} promoted the view that the
Hilbert space $\H$ on which a representation $(\pi,S)$ is defined
should be considered as a Hilbert module over $\rC_0(X) \times_\sigma \Fn$.
Muhly and Solel \cite{MSmem} have adopted this view,
although they changed the nomenclature somewhat.
They focussed on two notions:
A Hilbert module $\K$ over an operator algebra $\A$
is \textit{orthogonally injective} if every contractive short exact sequence
\[ 0 \to \K \to \M \to \Q \to 0 \]
has a contractive splitting.  Likewise a Hilbert module $\Q$
is \textit{orthogonally projective} if every such contractive
short exact sequence has a contractive splitting.
It is not difficult to show that these two properties together
are equivalent to being a maximal representation.
This means that Muhly and Solel's characterization \cite{MSbound} of
the unique extension property and (not necessarily irreducible)
boundary representations is essentially the same as that
of Dritschel--McCullough and Arveson.

If $\pi$ is a completely contractive representation of
$\rC_0(X) \times_\sigma \Fn$, then it is not difficult
to show that $\H_\pi$ is orthogonally projective if and only if
$\pi$ is isometric.  The details are left to the interested reader.

\begin{prop}\label{orthoinjective}
Let $(\pi,S)$ be a contractive representation of the semicrossed product
$\rC_0(X) \times_\sigma \Fn$ on a Hilbert space $\H$,
considered as a Hilbert $\rC_0(X) \times_\sigma \Fn$ module.
The following are equivalent:
\begin{enumerate}
\item $\pi$  is a full isometric representation.
\item $\pi$ is a maximal representation.
\item The $S_i$'s are isometries and $\H$ is orthogonally injective.
\end{enumerate}
\end{prop}

\begin{proof}
Suppose that (3) holds.
Let $(\rho,T)$ be  a full dilation of $\pi$ on a Hilbert space $\M$.
Since each $S_i$ is an isometry, one must have $T_i|\H = S_i$;
and hence $\H$ is invariant for $(\rho,T)$.
Therefore the complement $\Q = \H^\perp$ together with the restriction
$\tilde{\rho}$ of $\rho$ and the compressions $A_i = P_\Q T_i|_{\Q}$
is the quotient Hilbert module.  So
\[ 0 \to \H \to \M \to \Q \to 0 \]
is a contractive short exact sequence of Hilbert modules.
This splits, which implies that $\Q$ is invariant for each $S_i$.
Fullness of $(\rho,T)$ means that $T_i T_i^* = E_\rho(\sigma_i(X))$,
which means that
\[
 S_i S_i^* \oplus A_i A_i^* =
 E_\pi(\sigma_i(X)) \oplus E_{\tilde{\rho}}(\sigma_i(X)) .
\]
 From this it follows that each $S_i$ also satisfies
the fullness condition.  So (2) holds.

Now suppose that $(\pi,S)$ is a full isometric
representation on $\H$, and suppose that $(\rho,T)$
is a dilation of $(\pi,S)$ on $\M$.
As in the previous paragraph, the fact that each $S_i$ is an
isometry means that $\H$ is invariant.
By fullness, the range of each $S_i$ is all of
$E_\pi(\sigma_i(X)) \H$.  As the range of $T_i$ is contained
in $E_\rho(\sigma_i(X))\M$, it follows that $T_i \H^\perp$
is contained in $\H^\perp$.
Hence $\H$ reduces $(\rho,T)$; that is, $(\pi,S)$ is maximal.
So (1) holds.

The implication (1) implies (3) is trivial.
\end{proof}

It follows that the C*-envelope is obtained by taking all
fully isometric representations of $\rC_0(X) \times_\sigma \Fn$
on a Hilbert space of the appropriate size.

\begin{cor}\label{Cenv_semicross}
The C*-envelope of $\rC_0(X) \times_\sigma \Fn$ can be
obtained as $\ca(\rho(\A_0(X,\sigma)))$,
where $\rho$ is the direct sum of all full isometric representations
of $\A_0(X,\sigma)$ on a fixed Hilbert space of dimension
$\aleph_0 |X|$.  If $X$ is separable and metrizable, then
a separable Hilbert space will suffice.
\end{cor}

One useful conclusion is the following:

\begin{cor}\label{unitary elements}
If $\sigma_i$ is surjective, then $\fs_i$ is a unitary element
in $\cenv(\rC_0(X) \times_\sigma \Fn)$.
\end{cor}

\begin{eg}
It is not clear to us how to write down a complete family of
full isometric representations for the semicrossed product.
However,  the theorem above suggests that
one consider the \textit{full atomic representations} as a natural candidate.
For each $x\in X$, let $\H_x$ denote a non-zero Hilbert space such that
\[ \sum_{y \in \sigma_i^{-1}(x)} \dim \H_y = \dim \H_x \qforal x\in X \AND 1 \le i \le n .\]
Taking the dimension to be the
\[  \aleph_0 \max\{ |\sigma_i^{-1}(x)| : x\in X,\ 1 \le i \le n \}  \]
for every $\H_x$ will suffice.
Then for each $(x,i)\in X \times n$, select isometries $S_{i,y} \in \B(\H_y,\H_x)$
for each $y \in \sigma_i^{-1}(x)$ so that they have pairwise orthogonal ranges and
\[  \sotsum_{y \in \sigma_i^{-1}(x)} S_{i,y} S_{i,y}^* = I_{\H_x} .\]
Then define $T_i = \sotsum_{y \in X} S_{i,y}$.
Define $\pi(f)|_{\H_y} = f(y) I_{\H_y}$.
Then we have a full isometric representation of $(X,\sigma)$.

Interesting separable representations can be found by taking a dense
countable subset $Y$ such that $\sigma_i(Y) \subset Y$ and for every
$y \in Y \cap \sigma_i(X)$,  $Y \cap\sigma_i^{-1}(x)$ is non-empty.
Then build isometries as above.

Unfortunately, the general structure of such representations appears
to be very complicated.  This makes it difficult to describe the algebraic
structure in the C*-envelope.
\end{eg}

\begin{eg} \textbf{Case $n=1$ with surjective map.}
Suppose that $n=1$, $X$ is compact, and that $\sigma$ is surjective.
By Corollary~\ref{unitary elements}, the semicrossed
product $\rC(X) \times_\sigma \bZ^+$ is a subalgebra of
its C*-envelope generated by a copy of $\rC(X)$ and a
\textit{unitary} operator,  even if $\sigma$ is a multiple to one map.
In \cite{Pet2}, Peters constructs a canonical solenoid type space $Y$
with a projection $p$ of $Y$ onto $X$ and a homeomorphism
$\tau$ of $Y$ such $\sigma p = p \tau$.
He shows that the C*-envelope of $\rC(X) \times_\sigma \bZ^+$
is the crossed product C*-algebra $\rC(Y) \times_\tau \bZ$.
\end{eg}

Finally, we show that in general the algebra $\rC_0(X) \times_\sigma \Fn$ does not
arise as the tensor algebra of some C*-correspondence. This implies in
particular that the computations of Section \ref{S:correspondence} have no analogues
in the context of semicrossed products when $n >1$.

\begin{prop} \label{imposs}
Let $(X, \sigma)$, $\sigma=(\sigma_1, \sigma_2, \dots,
\sigma_n)$, $n>1$, be a multivariable dynamical system.
Assume that the maps $\sigma_i$ have a common fixed point.
Let $\E$ be an arbitrary C*-correspondence, over a C*-algebra $\fA$,
and let $\T_+(\E)$ be the associated tensor algebra.
Then, the semicrossed product $\rC_0(X) \times_\sigma \Fn$
and the tensor algebra $\T_+(\E)$ are not completely isometrically isomorphic.
\end{prop}

\begin{proof}
Let $x_0$ be a common fixed of the $\sigma_i$ for $1 \le i \le n$.
Consider the representation $\pi$ of $\rC_0(X) \times_\sigma \Fn$ which sends
$f \in \rC_0(X)$ to $f(x_0)$ and sends $\fs_i$ to the generators $U_i$ of $\ca_r(\bF_n)$,
the reduced C*-algebra of the free group $\bF_n$.
This representation is full since each $U_i$ is unitary.
Therefore by Proposition~\ref{orthoinjective}, $\pi$ is maximal.
Hence $\pi$ has a unique completely positive extension
to $\cenv(\rC_0(X) \times_\sigma \Fn)$ and it is a $*$-representation.
Thus $\ca_r(\bF_n)$ is a quotient of $\ca(\rC_0(X) \times_\sigma \Fn)$.
In particular, the C*-envelope is not nuclear.

On the other hand, suppose that $\rC_0(X) \times_\sigma \Fn$
and $\T_+(\E)$ were completely isometrically isomorphic.
The diagonal of $\T_+(\E)$, and hence the diagonal of
$\rC_0(X) \times_\sigma \Fn$, contains an isomorphic copy of $\fA$.
Therefore $\fA$ is commutative.
Consequently the associated Cuntz-Pimsner algebra $\O(\E)$ is
nuclear \cite[Corollary~7.5]{Ka}.
By a result of Katsoulis and Kribs \cite{KK2},
$\O(\E)$ is the C*-envelope of $\T_+(\E)$.

Thus as the C*-envelope of $\rC_0(X) \times_\sigma \Fn$ is not nuclear,
it is not a tensor algebra.
\end{proof}

For specific correspondences $\E$, one can draw stronger
conclusions than that of the above Proposition.
For instance we will see in Corollary~\ref{C:noniso} that under
weaker hypotheses than Proposition \ref{imposs}, the semicrossed product
$\rC_0(X) \times_\sigma \Fn$ and the tensor algebra $\A(X, \sigma)$ are not
isomorphic \textit{as algebras}.

\chapter{Recovering the Dynamics} \label{C:class}

\section{Fourier Series and Automatic Continuity}\label{S:fourier}

Since the tensor algebra and semicrossed product algebras
are defined by a universal property,  it is evident that
whenever $(\pi,S)$ satisfies the covariance relations and (row) contractivity,
then so does $(\pi, \lambda S)$ for $\lambda = (\lambda_i) \in \bT^n$,
where $\lambda S = \begin{bmatrix} \lambda_1 S_1 & \dots & \lambda_n S_n \end{bmatrix}$.
Therefore the map $\alpha_\lambda$ which sends the generators $\fs_i$
to $\lambda_i \fs_i$ and fixes $\rC_0(X)$ yields a completely isometric
isomorphism of $\A(X,\sigma)$ or $\rC_0(X) \times_\sigma \Fn$.
Consequently it extends uniquely to a $*$-automorphism of the C*-envelope.
In particular, if $\lambda_i = z \in \bT$ for $1 \le i \le n$, we obtain the
gauge automorphisms $\gamma_z$ that are a key tool in many related studies.

One immediate application is standard:

\begin{prop}\label{expectation}
The map $E(a) = \dint_\bT \gamma_z(a) \,dz$ \vspace{.2ex}
is a completely contractive expectation of $\A(X,\sigma)$
and $\rC_0(X) \times_\sigma \Fn$ onto $\rC_0(X)$.
\end{prop}

\begin{proof} It is routine to check that the map taking $z$ to $\gamma_z(a)$ is norm
continuous.  Check this on $\A_0(X,\sigma)$ first and then approximate.
So $E(a)$ makes sense as a Riemann integral.
Now computing $E$ on the monomials $\fs_w f$ shows that $E(f)=f$ and $E(\fs_w f)=0$
for $w \ne \mt$.
As this map is the average of completely isometric maps, it is completely contractive.
\end{proof}

The map $E$ as defined above makes sense for any element of
the C*-envelope.
However the range is then not $\rC_0(X)$, but rather the span of all words
of the form $\fs_v f \fs_w^*$ for $|v|=|w|$.
We will make use of this extension below.

We can see this explicitly for the tensor algebra.
Observe that for the boundary representations of Theorem~\ref{T:boundary},
one can see the expectation as the compression to the diagonal.
It is follows from this representation that one can read off
the \textit{Fourier coefficients} of elements of $\A(X,\sigma)$
This is computed within the C*-envelope as follows:

\begin{defn}
For each word $w \in \Fn$, define a map $E_w$ from
$\A(X,\sigma)$ onto $\rC_0(X)$ by $E_w(a) = E(\fs_w^* a)$.
\end{defn}

Observe that $\fs_w^* \fs_v = \fs_u$ when $v=wu$, it equals $\fs_u^*$ when $w = vu$,
and otherwise it is $0$.
Thus the only time that $E(\fs_w^* \fs_v f) \ne 0$ is when $v=w$.
Hence if a polynomial $a = \sum_{v\in\Fn} \fs_v f_v \in \A_0(X,\sigma)$,
it is clear that $E_w(a) = f_w$.

These Fourier coefficients do not seem to be obtainable from an integral
using invariants of $\A(X,\sigma)$ without passing to the C*-algebra.
This means that they are less accessible than in the singly generated case.
A partial recovery is the following.  For $k \ge 0$, define
\[ \Phi_k(a) = \int_\bT \gamma_z(a) \ol{z}^k \,dz .\]
This is clearly a completely contractive map.
Checking it on monomials shows that
\[ \Phi_k(\fs_w f) = \begin{cases} \fs_w f &\qif |w|=k\\ 0 &\qif |w| \ne k \end{cases} .\]
Indeed, one can do somewhat better and obtain a sum over
all words with the same abelianization.
That is, if $\bk = (k_1,\dots ,k_n) \in \bN_0^n$, define
\[ \Psi_\bk(a) = \int_{\bT^n} \alpha_\lambda(a) \ol{\lambda}^\bk \,d\lambda \]
where $ \ol{\lambda}^\bk = \prod_{i=1}^n \ol{\lambda}_i^{k_i}$.
Again it is easy to check that this is a completely contractive map
onto the span of words $\fs_w f$ such that $w(\lambda) = \lambda^\bk$.

As for Fourier series, this series generally does not converge
for arbitrary elements of $\A(X,\sigma)$.
However one can define the Cesaro means and recover $a$ from its Fourier series.
Define the $k$th Cesaro mean by
\[
 \Sigma_k(a) = \sum_{i=0}^k \big( 1 - \tfrac i k \big) \sum_{|w|=i} \fs_w E_w(a)
 =  \sum_{i=0}^k \big( 1 - \tfrac i k \big) \Phi_i(a)  .
\]
As usual, this may be obtained as an integral against the Fejer kernel $\sigma_k$ by
\[
 \Sigma_k(a) = \int_\bT \gamma_z(a) \sigma_k(z) \,dz .
\]
Since $\sigma_k$ is positive with $\|\sigma_k\|_1 = 1$, this is again
a completely contractive map.
A routine modification of the usual Fejer Theorem of classical
Fourier analysis shows that
\[ a = \lim_{k\to\infty} \Sigma_k(a) \qforal a \in \A(X,\sigma) . \]
So we may write $a \sim \sum_{w \in \Fn} \fs_w f_w$,
where $f_w = E_w(a)$, to mean that this is the Fourier series
of $a$, with summation interpreted via the Cesaro means.

A useful fact that derives from the Fourier series is the following:

\begin{prop}\label{P:factor}
Fix $k \ge 1$ and $a \in \A(X,\sigma)$.
Suppose that $E_v(a) = 0$ for all $|v| < k$.
Then $a$ factors as $a = \sum_{|w|=k} \fs_w a_w$
for elements $a_w \in \A(X,\sigma)$;
and $\|a\| = \big\|  \sum_{|w|=k} a_w^* a_w \big\|^{1/2}$.
\end{prop}

\begin{proof}  First suppose that $a \in \A_0(X,\sigma)$.
Observe that $a_w := \fs_w^* a$ belongs to $\A(X,\sigma)$.
It is then clear that
\[ \sum_{|w|=k} \fs_w a_w =\big( \sum_{|w|=k} \fs_w \fs_w^*\big) a = a .\]
Moreover $a$ factors as
\[ a = \begin{bmatrix} \fs_{w_1} & \dots & \fs_{w_{n^k}} \end{bmatrix}
 \begin{bmatrix}a_{w_1}\\ \vdots\\ a_{w_{n^k}}  \end{bmatrix}
\]
where $w_1,\dots,w_{n^k}$ is any enumeration of the words of length $k$.
Since $ \begin{bmatrix} \fs_{w_1} & \dots & \fs_{w_{n^k}} \end{bmatrix}$
is an isometry, if follows that
\[ \|a\| = \big\|  \sum_{|w|=k} a_w^* a_w \big\|^{1/2} . \]

For an arbitrary $a \in \A(X,\sigma)$ with $E_v(a) = 0$ for all $|v| < k$,
note that its Cesaro means have the same property.
Hence we can similarly define $a_w = \fs_w^* a$ and verify
the result by taking a limit using these polynomials.
\end{proof}

For an arbitrary element of $\A(X,\sigma)$, subtract off the Fourier
series up to level $k-1$ and apply the proposition.  One gets:

\begin{cor}\label{C:factor}
Fix $k \ge 1$ and $a \in \A(X,\sigma)$.
Then $a$ can be written as $a = \sum_{|v|<k} \fs_v E_v(a) + \sum_{|w|=k} \fs_w a_w$
for certain elements $a_w \in \A(X,\sigma)$.
\end{cor}

The isometries in $\rC_0(X) \times_\sigma \Fn$ do not have nice relations,
so one cannot define the Fourier coefficients as easily as in the
tensor algebra case.
One can, though, define the projections $\Phi_k$ onto the span of
all words of length $k$ in $T_1,\dots,T_n$.
This is a finite dimensional subspace spanned by these words,
which form a basis.  So in principle, there is a Fourier series
that can be determined in this way.
Certainly the Cesaro means exist as nice integrals just as before.
These are completely contractive maps into $\A_0(X,\sigma)$
which converge in norm for every element of $\rC_0(X) \times_\sigma \Fn$.
There is no analogue of Proposition~\ref{P:factor}.

An important feature of the Fourier series expansion
is that if $\Phi_k(A) = 0$ for all $k \ge0$, then $A=0$
in both $\A(X,\sigma)$ and $\rC_0(X) \times_\sigma \Fn$.
This is a key property in establishing automatic continuity for
isomorphisms..

Recall that if $\phi : \A \to \B$ is an epimorphism between
Banach algebras, then the \textit{separating space} of $\phi$ is the
two-sided closed ideal of  $\B$ defined as
\[
\S (\phi ) := \big\{ b \in \B  : \exists \{a_n\}_{n\ge1} \subseteq \A
\mbox{ so that } a_n \to 0 \mbox{ and } \phi(a_n )\to b  \big\}.
\]
Clearly the graph of $\phi$ is closed if and only if $\S (\phi )= \{0\}$.
Thus by the closed graph theorem, $\phi$ is continuous
if and only if $\S (\phi ) = \{0\}$.

The following is an adaption of \cite[Lemma 2.1]{Sinclair} and was used in
\cite{DKisotop}.

\begin{lem}[Sinclair]     \label{Sinclair}
Let  $\phi : \A \to \B$ be an epimorphism between
Banach algebras and let $ \{b_k\}_{k\in \bN}$
be any sequence in $\B$.
Then there exists $k_0 \in \bN$ so that for all $k \ge k_0$,
\[
 \ol{ b_1 b_2 \dots b_k \S (\phi ) } = \ol{ b_1 b_2 \dots b_{k+1}  \S (\phi ) }
\]
and
\[
 \ol{ \S (\phi ) b_{k} b_{k-1} \dots b_1 } = \ol{\S (\phi ) b_{k+1} b_k \dots b_1 }.
\]
\end{lem}

\begin{cor} \label{auto_cont}
Let $(X, \sigma)$ and $(Y, \tau)$ be multivariable dynamical systems.
Then any isomorphism $\gamma$ of $\A(X, \sigma)$ onto $\A(Y , \tau)$
or of $\rC_0(X) \times_\sigma \Fn$ onto $\rC_0(Y) \times_\tau \Fn$
is automatically continuous.
\end{cor}

\begin{proof}
Fix one of the generating isometries of $\A(Y , \tau)$, say $\ft_1$.
For any subset $\S$ of $\A(Y , \tau)$,
the faithfulness of the Fourier series expansion implies that
\[ \bigcap_{k\ge0} \ft_1^k \S  = \{0\} . \]
Thus if $\S(\gamma) \neq \{0\}$, then taking
$b_i =\ft_1$ in Lemma~\ref{Sinclair}, we obtain an integer $k_0$ so that
\[ \ft_1^{k_0}\S(\gamma)  = \bigcap_{k\ge0} \ft_1^k \S(\gamma) = \{0\} .\]
Since left multiplication by $\ft_1$ is injective, $\S(\gamma) = \{0\}$.
Therefore $\gamma$ is continuous.

The same argument works for the semicrossed product.
\end{proof}

This result allows us to consider only \textit{continuous} representations in the study
of arbitrary isomorphisms between tensor algebras of multisystems.

\section{Characters and Nest Representations.}\label{S:characters}

In this section, we extend methods from \cite{DKisotop} to the multivariable setting.
This will be applied in the Section~\ref{S:classification} to recover much
of the dynamical system from the tensor or crossed product algebra.
Following Hadwin-Hoover \cite{HH1}, we first look at characters.
At fixed points, there will be some analytic structure which will be important.
Then we study nest representations into the $2\times 2$ upper triangular matrices $\fT_2$.

In this section, many results will hold for both the tensor algebra $\A(X,\sigma)$
and the semicrossed product $\rC_0(X) \times_\sigma \Fn$.
We will use $\A$ to denote either algebra, and will specify when the results diverge.

\subsection*{Characters.}
Let $\fM_\A$ denote the space of characters of $\A$ endowed with the weak-$*$ topology.
Since $\A$ contains $\rC_0(X)$ as a subalgebra, the restriction of any character
$\theta$ to $\rC_0(X)$ will be a point evaluation $\delta_x$ at some point $x \in X$.
Let $\fM_{\A,x}$ denote the set of all characters extending $\delta_x$.
Observe that since there is an expectation $E$ of $\A$ onto  $\rC_0(X)$,
there is always a distinguished character $\theta_{x,0} = \delta_x E$ in $\fM_{\A,x}$.
Since characters are always continuous, a character $\theta \in \fM_{\A,x}$ is determined
by $z = (\theta(\fs_1),\dots,\theta(\fs_n))$.
We will write $\theta_{x,z}$ for this character when it is defined.

\begin{lem}\label{L:characters}
Let $x \in X$, and let $I_x = \{ i : 1 \le i \le n,\ \sigma_i(x) = x \}$.
Then
\[ \fM_{\A(X,\sigma),x} = \{ \theta_{x,z} : z_i = 0 \FOR i \notin I_x,\ \|z\|_2 \le 1 \} =: \ol{\bB(I_x)} \]
and
\[ \fM_{\rC_0(X) \times_\sigma \Fn,x} = \{ \theta_{x,z} : z_i = 0 \FOR i \notin I_x,\ \|z\|_\infty \le 1 \}
 =: \ol{\bD(I_x)} .\]
Moreover for each $a \in \A$, the function $\Theta_a(z) = \theta_{x,z}(a)$
is analytic on the ball (respectively polydisc) of radius 1 in the variables
$\{z_i : i \in I_x \}$ and is continuous on the closure.

In particular, $\fM_{\A,x} = \{ \theta_{x,0} \}$ if $x$ is not a fixed point for any $\sigma_i$.
\end{lem}

\begin{proof} Let $\theta \in \fM_{\A,x}$.
Characters always have norm 1; and indeed, they are completely contractive.
So for $\A(X,\sigma)$, we must have
\[
 \|z\|_2 = \| \theta^{(1,n)}(\fs) \| \le 1 .
\]
In case of the semicrossed product, we obtain $\|z\|_\infty \le 1$.

If $\sigma_i(x) = y \ne x$, then select a function $f\in\rC_0(X)$ such that
$f(x) = 0$ and $f(y) = 1$.  Then
\begin{align*}
 0 &= f(x)\theta(\fs_i) = \theta(f) \theta(\fs_i) = \theta(f\fs_i) \\
 &= \theta(\fs_i (f\circ\sigma_i)) = \theta(\fs_i) f(y) = \theta(\fs_i) .
\end{align*}

So now suppose that $z \in \ol{\bB(I_x)}$.
Define a one-dimensional representation of $\A(X,\sigma)$
by setting $\theta(f) = f(x)$ and $\theta(\fs_i) = z_i$.
Then the fact that $\sigma_i(x)=x$ for $i \in I_x$ ensures that
the covariance relations are satisfied.
Since $\|z\|_2 \le 1$, this is a row contractive representation.
Hence this extends to a contractive representation of $\A(X,\sigma)$,
yielding the desired character $\theta_{x,z}$.
Similarly if $z \in \ol{\bD(I_x)}$, then this determines
a contractive covariant representation.
So it extends to a character of $\rC_0(X) \times_\sigma \Fn$.

Now consider analyticity.  If $a \in \A_0(X,\sigma)$,
then $\Theta_a(z)$ is a polynomial in $\{z_i : i \in I_x \}$; and hence is analytic.
Since $\|\Theta_a\|_\infty \le \|a\|$, it now follows for arbitrary $a$
by approximation that $\Theta_a$ is the uniform limit of polynomials
on $\ol{\bB(I_x)}$ (respectively $\ol{\bD(I_x)}$).
Hence it is analytic on the interior of the ball (respectively polydisc)
and continuous on the closure.
\end{proof}

\begin{defn}
An open subset $M$ of $\fM_\A$ is called \textit{analytic} if there is a domain
$\Omega$ in $\bC^d$ and a continuous bijection $\Theta$ of $\Omega$ onto $M$
so that the function $\Theta(z)(a)$ is analytic on $\Omega$ for every $a \in \A$.
It is a \textit{maximal analytic subset} if it is maximal among analytic subsets of $\fM_\A$.
\end{defn}

In particular, consider $x \in X$ for which $I_x$ is non-empty.
The open ball $\bB(I_x) = \{ z \in \ol{\bB(I_x)} : \|z\|_2 < 1 \}$
considered as an open subset of $\bC^{|I_x|}$
is a complex domain that is mapped in the obvious way onto
$B_x = \{ \theta_{x,z} : z \in \bB(I_x) \}$.
This is an analytic set in $\fM_{\A(X,\sigma)}$.
Similarly the polydisc $\bD(I_x) = \{ z \in \ol{\bD(I_x)} : \|z\|_\infty < 1 \}$
maps onto the analytic set $D_x = \{ \theta_{x,z} : z\in \bD(I_x) \}$
in $\fM_{\rC_0(X) \times_\sigma \Fn}$.

\begin{lem}
The maximal analytic sets in $\fM_{A( X , \sigma)}$ are precisely the
balls $B_x$ for those points $x \in X$ fixed by at least one $\sigma_i$.
Similarly the maximal analytic sets in $\fM_{\rC_0(X) \times_\sigma \Fn}$ are precisely the
polydiscs $D_x$ for those points $x \in X$ fixed by at least one $\sigma_i$.
\end{lem}

\begin{proof} It suffices to show that analytic sets must sit inside one of the fibres $\fM_{\A,x}$.
Let $\Theta$ map a domain $\Omega$ into $\fM_\A$.
For every $f \in \rC_0(X)$, the function $\Theta(z)(f)$ and the function
$\Theta(z)(\ol{f}) = \ol{\Theta(z)(f)}$ are analytic, and hence constant.
Since continuous functions on $X$ separate points, this implies
that $\Theta$ maps into a single fibre.
 
Now observe that $\fM_{\A,x}$ is homeomorphic to a closed ball (respectively polydisc).
Hence the maximal analytic subset would be the ``interior'' $B_x$ (respectively $D_x$).
So we have identified all of the maximal analytic sets.
\end{proof}

\begin{cor}\label{C:quotient}
The characters of $\A$ determine $X$ up to homeomorphism, and identify
which points are fixed by some $\sigma_i$'s, and
determine exactly how many of the maps fix the point.
\end{cor}

\begin{proof} The lemma shows that $\fM_\A$ consists of a space which is fibred over $X$,
and the fibres are determined canonically as the closures of maximal
analytic sets and the remaining singletons.
Thus there is a canonical quotient map of $\fM_\A$ onto $X$;
and this determines $X$.
Next, the points which are fixed by some $\sigma_i$ are exactly the
points with a non-trivial fibre of characters.
The corresponding maximal analytic set is homeomorphic to a ball
(respectively polydisc) in $\bC^d$ where $d = |I_x|$.
The invariance of domain theorem shows that the dimension $d$
is determined by the topology.
\end{proof}

\begin{cor}\label{C:noniso}
If $(X,\sigma)$ has a point fixed by two or more of the maps $\sigma_i$,
then $\A(X,\sigma)$ and $\rC_0(X) \times_\sigma \Fn$ are not
algebraically isomorphic.
\end{cor}

\begin{proof}
An algebra isomorphism will yield a homeomorphism of the
character spaces, and will be a biholomorphic map of each
maximal analytic set of the tensor algebra to the corresponding
maximal analytic set in the semicrossed product.
The existence of a point $x_0$ fixed by $k\ge2$ of the maps
$\sigma_i$ means that the tensor algebra contains the ball $\bB_k$
as a maximal analytic set, while the semicrossed product has
a polydisc $\bD_k$.
However no polydisk of dimension at least 2 is biholomorphic
to any ball, and vice versa.
Therefore the algebras are not isomorphic.
\end{proof}

\subsection*{Nest representations.}
In \cite{DKisotop}, we considered representations \textit{onto} the
$2\times2$ upper triangular matrices $\fT_2$.
Here we actually need to consider a more general notion.

\begin{defn}
Let $\N_2$ denote the maximal nest $\big\{ \{0\}, \bC e_1, \bC^2 \big\}$
in $\bC^2$.
If $\A$ is an operator algebra, let $\rep_{\N_2}$ denote the
collection of all continuous representations  $\rho$ of $\A$ on $\bC^2$
such that $\Lat \rho(\A) = \N_2$.
\end{defn}

These representations are called \textit{nest representations}.
There are two unital subalgebras of the $2\times2$ matrices $\fM_2$
with this lattice of invariant subspaces, $\fT_2$ and the
abelian algebra $\A(E_{12}) = \spn\{ I, E_{12} \}$.
Both have non-trivial radical.
The other unital subalgebras of $\fT_2$ are semisimple, and except for $\bC I$,
they are  all similar to the diagonal algebra $\fD_2$.
Their lattice of invariant subspaces is the Boolean algebra with two generators.
So representations with semisimple range are not nest representations.

Observe that for any representation $\rho$ of $\A$ into $\fT_2$, the
compression to a diagonal entry is a homomorphism.
Thus $\rho$ determines two characters which we denote by
$\theta_{\rho,1}$ and $\theta_{\rho,2}$.
The map $\psi$ taking $a$ to the $1,2$-entry of $\rho(a)$
is a point derivation satisfying
\[
 \psi(ab) = \theta_{\rho,1}(a) \psi(b) + \psi(a) \theta_{\rho,2}(b) \qfor a,b \in \A .
\]
Define $\rep_{y,x}(\A)$ to be those nest representations
for which $\theta_{\rho,1} \in \fM_{\A,y}$ and $\theta_{\rho,2} \in \fM_{\A,x}$.

It is convenient to consider representations which restrict
to $*$-repre\-sent\-ations on $\rC_0(X)$.
For (completely) contractive representations of $\rC_0(X)$,
this is automatic.
It is also the case that representations of $\rC_0(X)$
into $\fM_2$ are automatically continuous, and thus
are diagonalizable. Here we need a stronger version of this fact.
Let $\rep^d \A$ and $\rep^d_{y,x} \A$ denote the nest representations
which are diagonal on $\rC_0(X)$.

\begin{lem} \label{diagonalizable}
Let $X$ be a locally compact space; and let $\sigma$ be a continuous
map of $X$ into itself.
Let $K \subset X$, and let $\Omega$ be a domain in $\bC^d$.
Suppose that $\rho$ is a map from $K \times \Omega$ into
$\rep\A$ so that
\begin{enumerate}
\item $\rho_{x,z} := \rho(x,z) \in \rep_{x,\sigma(x)}\A$ for $x \in K$ and $z \in \Omega$.
\item $\rho$  is continuous in the point--norm topology.
\item For each fixed $x \in K$, $ \rho(x,z)$ is analytic in $z \in \Omega$.
\end{enumerate}
Then there exists a map $A$ of $K \times \Omega$ into $\fT_2^{-1}$,
the group of invertible upper triangular matrices,
so that
\begin{enumerate}
\item $A(x, z)\rho_{x,z}(\cdot)A(x, z)^{-1} \in \rep^d \A$.
\item $A(x,z)$ is continuous on $\big( K\bsl \{x: \sigma(x)=x \} \big) \times \Omega$.
\item For each fixed $x \in K$, $A(x,z)$ is analytic in $z \in \Omega$.
\item $\max\{ \|A(x,z)\|, \|A(x,z)^{-1} \| \} \le 1 + \| \rho_{x,z} \|$.
\end{enumerate}
\end{lem}

\begin{proof}
When $\sigma(x)=x$, every representation $\rho \in \rep_{x,x}\A$
satisfies $\rho(f) = f(x) I_2$, and so they are automatically diagonalized.
Define $A(x,z) = I_2$. This clearly satisfies conclusions (1), (3) and (4).

Now consider a point $x$ with $\sigma(x) \ne x$.
Choose a compact neighbourhood $\ol{\V}$ of $x$
so that $\sigma(\ol{\V})$ is disjoint from $\ol{\V}$.
Select $f \in \rC_0(X)$ so that $f(\sigma(\ol{\V}))=\{1\}$ and $f(\ol{\V})=\{0\}$.
Let $F$ be the function on $\ol{\V} \times \Omega$ given
by the $1,2$ entry of $\rho_{x,z}(f)$; i.e.
\[
\rho_{x,z}(f)=\begin{bmatrix}1&F(x,z)\\0&0 \end{bmatrix}.
\]
Clearly $F$ is continuous, and is analytic in $z$ for each point in $\ol{\V}$.

The function $F$ does not depend on the choice of $f$; for if $g$
is another such function, the difference $f-g$ vanishes on
$\ol{\V} \cup \sigma(\ol{\V})$.
It is routine to factor this as the product $h_1h_2$ of two functions
vanishing on this set.
But then
\[ \rho_{x,z}(f-g) = \rho_{x,z}(h_1) \rho_{x,z}(h_2) \]
is the product of two strictly upper triangular matrices, and hence is $0$.
In particular, we may suppose that $\|f\| = 1$.

We can now define a function on
$\big( K\bsl \{x: \sigma(x)=x \} \big) \times \Omega$ by
\[
A(x,z) := \begin{bmatrix}1&F(x,z)\\0&1 \end{bmatrix} .
\]
This is well defined by the previous paragraph.
It is continuous, and is analytic in the second variable,
because it inherits this from $F$.

For a point $x$ with $\sigma(x) \ne x$, use the function $f$
selected above and compute
\begin{align*}
 A(x ,z)\rho_{x,z}(f)A(x,z)^{-1} &=
 \begin{bmatrix}1&F(x,z)\\0&1 \end{bmatrix}
 \begin{bmatrix}1&F(x,z)\\0&0 \end{bmatrix}
 \begin{bmatrix}1&-F(x,z)\\0&1 \end{bmatrix} \\&=
 \begin{bmatrix}1& 0\\0&0 \end{bmatrix} .
\end{align*}
Similarly, if $h\in\rC_0(X)$ and $h(\ol{\V}\, \cup\, \sigma(\ol{\V})) = \{1\}$,
then
\[ A(x,z)\rho_{x,z}(h)A(x,z)^{-1}= I_2 .\]
Since the image of $\rC_0(X)$ under any representation in $\rep_{x,\sigma(x)}\A$
is at most two dimensional, we conclude that
$A(x,z)\rho_{x,z}(.)A(x,z)^{-1}$ is diagonal when restricted to $\rC_0(X)$.

Finally observe that
\[ \|A(x,z) \| = \|A(x,z)^{-1}\| \le 1 +  \| \rho_{x,z} \| \|f\| = 1 + \| \rho_{x,z} \| . \qedhere \]
\end{proof}

\bigskip
\begin{eg}
For each $x \in X$ and $1 \le i \le n$, we define a
nest representation $\rho_{x, i}$ by
\[
 \rho_{x,i}(f) =  \begin{bmatrix}  f(\sigma_i(x)) &0\\  0 &f(x) \end{bmatrix}
 \qand \rho_{x,i}(\fs_j) = \begin{bmatrix}  0 &\delta_{ij} \\  0 & 0 \end{bmatrix}
\]
where $\delta_{ij}$ is the Kronecker delta function.  Thus
\[
 \rho_{x,i} \big( \sum_{w \in \Fn} \fs_w f_w  \big) =
 \begin{bmatrix}
  f_\mt(\sigma_i(x)) & f_i(x)\\
  0 & f_\mt(x)
 \end{bmatrix}
\]
This has $\theta_{\rho_{x,i},1} = \theta_{\sigma_i(x),0}$
and $\theta_{\rho_{x,i},2} = \theta_{x,0}$.
The $1,2$-entry $\psi_{\rho_{x,i}}(a) = f_i(x)$ is easily
see to satisfy the derivation condition.
This yields a nest representation of $\A$.
Since one generator is sent to a contraction and the rest
are sent to $0$, this clearly defines a completely contractive
representation of both $\A(X,\sigma)$ and $\rC_0(X) \times_\sigma \Fn$.

This representation maps $\A$ onto $\fT_2$ if $\sigma_i(x) \ne x$,
and onto $\A(E_{12})$ when $\sigma_i(x) = x$.
\end{eg}

The key to recovering the system $(X,\sigma)$ from $\A$ is the following lemma.

\begin{lem}\label{L:xtoy}
If $\rep_{y,x}(\A)$ is non-empty, then there is some
$i$ so that $\sigma_i(x) = y$. Furthermore, if $\rho \in \rep^d_{y,x}(\A)$
and $\sigma_j(x) \neq y$, then $\rho(\fs_j g)$ is diagonal for all $g \in \rC_0(X)$.
\end{lem}

\begin{proof}
Let $\rho \in \rep_{y,x}(\A)$.
By Lemma~\ref{diagonalizable}, we may assume
that $\rho$ belongs to $\rep^d_{y,x}(\A)$.
So for $f \in \rC_0(X)$, we have
$\rho(f) =  \begin{sbmatrix}  f(y) &0\\  0 &f(x) \end{sbmatrix}$.
If $\sigma_j(x) \ne y$, choose a function $f \in \rC_0(X)$
so that $f(y)=1$ and $f(\sigma_j(x)) = 0$.
Let $g \in \rC_0(X)$ and $s$ be the $1,2$ entry of $\rho(\fs_j g)$.
Apply $\rho$ to the identity $f \fs_j g= \fs_j g (f\circ\sigma_j)$
to obtain
\begin{align*}
\begin{bmatrix}*&s\\0&* \end{bmatrix} &=
\begin{bmatrix}f(y)\!&0\\0&\!f(x) \end{bmatrix}
\begin{bmatrix}*&s\\0&* \end{bmatrix} \\ &=
\begin{bmatrix}*&s\\0&* \end{bmatrix}
\begin{bmatrix}f(\sigma_i(y))\!\!&0\\0&\!\!f(\sigma_i(x)) \end{bmatrix} =
\begin{bmatrix}*&0\\0&0\end{bmatrix}
\end{align*}
Hence $s=0$, and so $\rho(\fs_j g)$ is diagonal for all $g \in \rC_0(X)$.

If $\rho(\fs_i g)$, $g \in \rC_0(X)$, are all diagonal for $1 \le i \le n$, then
$\rho(\A)$ is diagonal---and thus not a nest representation.
So there is some $i$ for which $\rho(\fs_i g)$ is not diagonal;
and thus $\sigma_i(x)=y$.
\end{proof}

\section{Piecewise conjugate multisystems}\label{S:conjugacy}

In this section we introduce the concept of piecewise conjugacy
for multivariable dynamical systems. This concept is central for
the classification of their algebras and apparently new in the theory of
dynamical systems. We therefore study some of its basic properties
and present a few illuminating examples.

Recall that two multivariable dynamical systems $(X, \sigma)$
and $(Y, \tau)$ are said to be \textit{conjugate} if there exists a
homeomorphism $\gamma$ of $X$ onto $Y$
and a permutation $\alpha \in S_n$ so that
$\gamma^{-1} \tau_i \gamma =  \sigma_{\alpha(i)}$ for $1 \le i \le n$.

\begin{defn}
Say that two multivariable dynamical systems $(X, \sigma)$ and $(Y, \tau)$ are
 \textit{piecewise conjugate} if there is a
homeomorphism $\gamma$ of $X$ onto $Y$
and an open cover $\{ \U_\alpha : \alpha \in S_n \}$ of $X$
so that for each $\alpha \in S_n$,
\[
 \gamma^{-1} \tau_i \gamma|_{\U_\alpha} = \sigma_{\alpha(i)} |_{\U_\alpha} .
\]
\end{defn}

We could have expressed this somewhat differently, saying that $Y$ has
an open cover so that on each open set, there is some permutation $\alpha$ so that
$\gamma$ locally intertwines each $\tau_i$ with $ \sigma_{\alpha(i)}$.
But this is readily seen to be equivalent.

The difference in the two concepts of conjugacy lies on the fact that
the permutations depend on the particular open set.
As we shall see, a single permutation generally will not suffice.

For a point $x \in X$, we say that two continuous functions $f,g$ mapping
$X$ into another topological space $Z$ are equivalent if they agree on an
open neighbourhood of $x$.
The equivalence class $[f]_x$ is called the \textit{germ} of $f$ at $x$.

Let $\sigma(x) = \{ \sigma_1(x),\dots.\sigma_n(x) \}$.

\begin{prop} \label{P:connected}
Let $(X, \sigma)$ and $(Y, \tau)$ be piecewise conjugate
multivariable dynamical systems.
Assume that $X$ is connected and that
 $E:= \{ x \in X : |\sigma(x)|=n \}$ is dense in $X$.
Then $(X, \sigma)$ and $(Y, \tau)$ are conjugate.
\end{prop}

\begin{proof}
Let $\gamma$ be the homeomorphism implementing the piecewise conjugacy.
Observe that by continuity, $E$ is open; and by hypothesis, it is dense.
For each $\alpha \in S_n$, define an open set
\[
 E_\alpha =
 \{ x \in X : [\gamma^{-1} \tau_i \gamma]_x = [\sigma_{\alpha(i)}]_x  \FOR 1 \le i \le n \} .
\]
 From the piecewise conjugacy, this is an open cover of $X$.
If $x\in E$, the permutation $\alpha$ must be unique since $\sigma(x)$
consists of $n$ distinct points.
This must persist to the closure, for if $E_\alpha \cap E_\beta$
is non-empty, then being open, it would intersect $E$,
contrary to the previous conclusion.
Hence the sets $E_\alpha$ must be clopen.
From the connectedness of $X$, we conclude that
there is a single permutation $\alpha$ so that $X=E_\alpha$.
This yields the desired conjugacy.
\end{proof}

Let $\Z(X, \sigma) = \{ x \in X : |\sigma(x)| < n \}$.
Then $\Z(X,\sigma)$ is the closed set $E^c$.
Thus the previous proposition could be stated
as saying that $\Z(X, \sigma)$ is nowhere dense.

If $\Z(X, \sigma)$ has non-empty interior,
then the situation may change dramatically.
We illustrate this in the simplest possible case.

\begin{eg}
Let $X=[0,1]$ and $\sigma= (\sigma_1, \sigma_2)$.
We describe the piecewise conjugacy class of $([0,1], \sigma)$.
We may as well assume that $Y=X$ and $\gamma = \id$.

Pick points $a_k , b_k$, $k \in K$, so that the interior of $\Z(X, \sigma)$
can be expressed as a \textit{disjoint} union
\[
 \Z(X, \sigma)^{\circ} = \bigcup_{k \in K} (a_k , b_k)
\]
of (relatively) open subintervals of $[0,1]$.
Select finitely many indices $k_1, k_2, \dots, k_l$.
Without loss of generality,
assume that $b_{k_j}< a_{k_{j+1}}$ for all $j=1, 2,\dots, l-1$.
For convenience, set $a_0=b_0=0$ and $a_{l+1} =b_{l+1}= 1$.
Define
\begin{align*}
 \tau_1(x) &= \begin{cases} \sigma_1(x) &\qif a_{2j} \le x \le b_{2j+1}\\
                      \sigma_2(x) &\qif a_{2j-1} \le x \le b_{2j}\end{cases}
\\ \intertext{and}
 \tau_2(x) &= \begin{cases} \sigma_2(x) &\qif a_{2j} \le x \le b_{2j+1}\\
                      \sigma_1(x) &\qif a_{2j-1} \le x \le b_{2j}\end{cases}
\end{align*}
Then $\U_\id = \bigcup (a_{2j}, b_{2j+1})$ and $\U_{(12)} = \bigcup (a_{2j-1}, b_{2j})$
is the partition for the piecewise conjugacy.

Conversely, any system $([0,1], \tau)$ which is piecewise conjugate to
$([0,1], \sigma)$ arises that way. The main point is that the functions
must coincide on an interval in order to make the switchover.
On the other hand, there cannot be countably many switches, because
the switching points would then have a cluster point, and  the two systems
could not coincide on any neighbourhood of that point.
We omit the details.

In particular, suppose that $\sigma_1$ is non-increasing, $\sigma_2$
is non-decreasing and that $\Z([0,1], \sigma)=[a,b]$, with $0<a<b<1$.
The switching across $(a,b)$ yields a piecewise conjugate system $(X,\tau)$
which is not monotonic.
As any homeomorphism of $[0,1]$ is monotone, these two systems are
definitely not conjugate.
\end{eg}

For $n=2$, we can be more definitive in Proposition~\ref{P:connected}.

\begin{prop}\label{P:connected_converse}
Let $X$ be connected and let $\sigma=(\sigma_1,\sigma_2)$;
and let $E:= \{ x \in X :|\sigma(x)|=2 \}$.
Then piecewise conjugacy coincides with conjugacy
if and only if $\ol{E}$ is connected.
\end{prop}

\begin{proof}
If $E$ is empty, then $\sigma_1=\sigma_2$ and there is nothing to prove.
As in the proof of Proposition~\ref{P:connected}, the sets $E_\alpha$
are open and
\[
 E_\id \cap E_{(12)} \cap \ol{E} = \mt .
\]
Arguing as before using the connectedness of $\ol{E}$,
it follows that there is some $\alpha$ so that
$E_\alpha$ contains $\ol{E}$.
But clearly both $E_\id$ and $E_{(12)}$ contain $\ol{E}^c = \Z(X,\sigma)^\circ$.
Thus $E_\alpha = X$ and so $X$ and $Y$ are conjugate.
\end{proof}

At the other end of the spectrum, total disconnectedness makes piecewise
conjugacy very tractable.

\begin{prop}\label{P:total_disc}
Let $X$ be a totally disconnected compact Hausdorff space.
Fix a homeomorphism $\gamma$ of $X$ onto another space $Y$.
Then there is a piecewise topological conjugacy of $(X,\sigma)$
and $(Y,\tau)$ implemented by $\gamma$ if and only if there
is a partition of $X$ into clopen sets $\{ \V_\alpha : \alpha\in S_n \}$
so that  for each $\alpha \in S_n$,
\[
 \gamma^{-1} \tau_i \gamma|_{\V_\alpha} = \sigma_{\alpha(i)} |_{\V_\alpha} .
\]
\end{prop}

\begin{proof}
If $(X,\sigma)$ and $(Y,\tau)$ are piecewise conjugate,
let
\[
 \U_\alpha =
  \{ x \in X : [ \gamma^{-1} \tau_i \gamma]_x = [\sigma_{\alpha(i)}]_x ,\ i \le i \le n \} .
\]
By hypothesis, this is an open cover of $X$.
Hence there is a partition of unity $\{ f_\alpha\}$ of positive
functions in $\rC(X)$ such that $\sum f_\alpha = 1$ and
$f_\alpha^{-1}(0,\infty) \subset \U_\alpha$ for $\alpha \in S_n$.
Then $K_\alpha = f_\alpha^{-1}[1/n!,1]$ is a compact subset
of $\U_\alpha$ and $\bigcup_\alpha K_\alpha = X$.

Each point $x \in K_\alpha$ is contained in $\U_\alpha$,
and thus there is a clopen neighbourhood $V_x$ of $x$
contained in $\U_\alpha$.
Select a finite subcover of $K_\alpha$.
The union of this finite cover is a clopen set $V_\alpha$
containing $K_\alpha$ and contained in $U_\alpha$.

To obtain the desired partition, order the permutations,
and replace $V_\alpha$ by its intersection with the complement of
the $V_\beta$'s which precede it in the list.
\end{proof}

\section{The Main Theorem}\label{S:classification}

The purpose of this section is to establish the connection between isomorphism
of the tensor or semicrossed product algebras and piecewise conjugacy.

Before embarking on the proof,  we need to recall some basic facts from
the theory of several complex variables \cite{FrG}.
Let $G \subseteq \bC^k$ be a domain and let $A\subseteq G$ be
an analytic variety of $G$.
One says that $A$ is regular at $z \in A$ of dimension $q$
if there is a neighborhood $U_z \subseteq G$ of $z$ and
holomorphic functions $f_i : U_z \to \bC$ for $1 \le i \le k-q$ so that
\begin{itemize}
\item[(i)] the Jacobian of $(f_1, f_2, \dots , f_{k-q})$ at $z$ has rank $k-q$, and,
\item[(ii)] $A \cap U_z = N(f_1, f_2, \dots f_{k-q}) := \bigcap_{i=1}^{k-q} f_i^{-1}(0).$
\end{itemize}
By the Theorem of local parameterization
at regular points \cite[Theorem I.8.3]{FrG}, the dimension is well defined.
The set $A_{\text{reg}}$ of all regular points of $A$ forms
an open subset of $A$, and the dimension is constant
on each connected component of $A_{\text{reg}}$.
The closure  of each such component is itself an
analytic variety of $G$ which is irreducible,
i.e., the regular points are connected.
The maximum dimension of the irreducible subvarieties of $A$
is called the dimension of $A$.

We need a rather modest conclusion from this theory.
However we know of no elementary proof of this fact.

\begin{prop} \label{zeroset}
Let $G \subseteq \bC^k$ be a domain and let
\[
 \psi =(\psi_1, \psi_2, \dots , \psi_l): \bC^k \to \bC^l, \text{  where  } l<k,
\]
be a holomorphic function.
Then the zero set $N(\psi_1, \psi_2, \dots , \psi_l)$
of $\psi$ is either empty or infinite.
In particular, the zero set has no isolated points.
\end{prop}

\begin{proof} Assume that $N(\psi_1, \psi_2, \dots , \psi_l)$ is non-empty.
If $l=1$, then the techniques of \cite[Proposition I.8.4]{FrG} imply
that the dimension of that variety is $k-1$.
Consequently, \cite[Exercise III.6.1]{FrG} implies that for an arbitrary $k$,
the dimension of $N( \psi_1, \psi_2, \dots , \psi_l)$ is at least $k-l$
and hence is positive.
In particular, some irreducible component of $N(\psi_1, \psi_2, \dots , \psi_l)$
is regular of positive dimension at some point $z$.
The conclusion follows now from the Implicit Function Theorem.

If there were an isolated point of the zero set, simply reduce $G$ to a
small neighbourhood of that point to obtain a contradiction.
\end{proof}

We now have all the requisite tools to recover much of the dynamical system
from the tensor or crossed product algebra in the  multivariable setting.
Again we will work as much as possible with both the tensor algebra
and the semicrossed product.

\begin{thm} \label{class}
Let $(X,\sigma)$ and $(Y,\tau)$ be two multivariable dynamical systems.
We write $\A$ and $\B$ to mean either
$\A(X,\sigma)$ and $\A(Y,\tau)$ or
$\rC_0(X)\times_\sigma\Fn$ and $\rC_0(Y)\times_\tau\Fn$.
If $\A$ and $\B$ are isomorphic as algebras, then the dynamical systems
$(X,\sigma)$ and $(Y,\tau)$ are piecewise conjugate.
\end{thm}


\begin{proof} Let $\gamma$ be an isomorphism of $\A$ onto $\B$.
This induces a bijection $\gamma_c$ from the character space
$\fM_\A$ onto $\fM_\B$ by $\gamma_c(\theta) = \theta\circ\gamma^{-1}$.
Similarly it induces a map $\gamma_r$ from $\rep_{\N_2}(\A)$
onto $\rep_{\N_2}(\B)$.

Since $\fM_\A$ is endowed with the weak-$*$ topology, it is
easy to see that $\gamma_c$ is continuous.
Indeed, if $\theta_\alpha$ is a net in $\fM_\A$ converging to $\theta$
and $b \in \B$, then
\[
 \lim_\alpha \gamma_c \theta_\alpha (b) = \lim_\alpha \theta_\alpha(\gamma^{-1}(a))
 = \theta(\gamma^{-1}(a)) = \gamma_c \theta(b) .
\]
The same holds for $\gamma_c^{-1}$.
So $\gamma_c$ is a homeomorphism.

Observe that $\gamma_c$ carries analytic sets to analytic sets.
Indeed, if $\Theta$ is an analytic function of a domain $\Omega$
into $\fM_\A$, then
\[ \gamma_c \Theta(z)(b) = \Theta(z)(\gamma^{-1}(b)) \]
is analytic for every $b \in \B$; and thus $\gamma_c \Theta$ is analytic.
Since the same holds for $\gamma^{-1}$, it follows that $\gamma_c$
takes maximal analytic sets to maximal analytic sets.
Thus it carries their closures, $\fM_{\A,x}$, onto sets $\fM_{\B,y}$.
The same also holds when these sets are singletons.

By Corollary~\ref{C:quotient}, $X$ is the quotient of $\fM_\A$
obtained by squashing each $\fM_{\A,x}$ to a point.
It follows that $\gamma_c$ induces a set map $\gamma_s$
of $X$ onto $Y$ which is a homeomorphism since
both $X$ and $Y$ inherit the quotient topology.

Fix $x _0 \in X$, and let $y_0= \gamma_c(x_0)$.
Fix one of the maps $\sigma_{i_0}$, and consider the set
\[
 \F =  \{ \sigma_i,\gamma^{-1} \tau_j \gamma
 : [\sigma_i]_x = [\sigma_{i_0}]_x =  [ \gamma^{-1} \tau_j \gamma]_x \} .
\]
For convenience, let us relabel so that $i_0=1$
and
\[
 \F = \{ \sigma_1,\dots,\sigma_k,
 \gamma^{-1} \tau_1 \gamma,\dots,\gamma^{-1} \tau_l \gamma\} .
\]
Fix a neighbourhood $\V$ of $x_0$ on which all of these
functions agree, and such that $\ol{\V}$ is compact. Furthermore,
if $\sigma_1(x_0)\neq x_0$, then choose $\V$ so that
$\ol{V} \cap \sigma_1(\ol{\V}) = \emptyset$.

The hard part of the proof is to show that $k=l$.
Assume that this has been verified.
Then one can partition the functions $\sigma_i$ into families with
a common germ at $x_0$, and they will be paired with a corresponding
partition of the $\gamma^{-1} \tau_j \gamma$'s of equal size.
This provides the desired permutation in some neighbourhood of each $x_0$.

It remains to show that $k=l$.
By way of contradiction, assume that $k \neq l$.
By exchanging the roles of $X$ and $Y$ if necessary,
we may assume that $k > l$.
We do not exclude the possibility that $l=0$.

For any $x \in \V$ and $z = (z_1, z_2, \dots , z_l) \in \bC^k$, consider
the covariant representations $\rho_{x, z}$ of $\A_0(X,\sigma)$ into $\fM_2$
defined by
\begin{align*}
\rho_{x,z}(f)&=
\begin{bmatrix}
  f(\sigma_1(x)) & 0\\
   0 &f(x)
\end{bmatrix},\\
\rho_{x,z}(\fs_i) &=
\begin{bmatrix}
  0 & z_i\\
   0 & 0
\end{bmatrix} \qfor 1 \le i \le k,\\
\intertext{and}
\rho_{x, z}(\fs_i) &=
\begin{bmatrix}
  0 & 0\\
  0 & 0
\end{bmatrix} \qfor k < i \le n.
\end{align*}

This extends to a well defined representation of $\A$, where a typical element
$A \sim \sum_{w \in \Fn} \fs_w f_w$ is sent to
\[
 \rho_{x,z}(A) =
 \begin{bmatrix}
  f_0(\sigma_1(x))  & \sum_{i=1}^k f_i(x) z_i\\
  0 & f_0(x)
 \end{bmatrix}
\]
There are no continuity problems since the Fourier coefficients are
continuous.

This representation will be (completely) contractive
if $z \in \ol{\bB}_k$ for $\A(X,\sigma)$ and $z \in \ol{\bD}_k$ for the semicrossed product.
For other values of $z$, these representations are similar
to completely contractive representations by conjugating by
$\diag(\|z\|_2,1)$ or $\diag(\|z\|_\infty,1)$ respectively.
Thus the norm can be estimated as $\|\rho_{x,z}\| \le \|z\|$,
where we use the 2-norm or the max norm depending on
whether we are considering the tensor algebra or the semicrossed product.

The representation $\rho_{x,z}$ maps into $\fT_2$ and is a nest representation
in $\rep_{\sigma_1(x), x} \A$ when $z \ne 0$, but is diagonal at $z=0$.
Observe that the range of $\rho_{x,z}$ for $z \ne 0$ equals $\fT_2$
when $\sigma_1(x) \ne x$ and equals $\A(E_{12})$ when $\sigma_1(x)= x$.
Moreover this map is point--norm continuous,
and is analytic in the second variable.

Now consider the map defined on $\ol{\V} \times \bC^k$ given by
\[
\Phi_0(x,z) = \gamma_r(\rho_{x,z})
 \in \rep_{\gamma_s\sigma_1(x), \gamma_s(x)}\B .
\]
By Corollary~\ref{auto_cont},
$\gamma$ is continuous;  and so $\gamma_r$ is also continuous.
Thus $\Phi_0$ is point--norm continuous,
and is analytic in the second variable.
So $\Phi_0$ fulfils the requirements of Lemma \ref{diagonalizable}.
Hence there exists a map $A(x, z)$
of $\ol{\V} \times \bC^k$ into $\fT_2^{-1}$, which
is analytic in the second variable, so that
\[
 \Phi(x, z) = A(x, z) \gamma_r(\rho_{x,z})  A(x, z)^{-1}
\]
diagonalizes $\rC_0(Y)$.
Moreover
\[ \max\{ \|A(x,z) \|, \| A(x,z)^{-1} \| \} \le 1 + \|\gamma_r\| \|z\| . \]
Recall that when $\sigma_1(x_0) \ne x_0$, we chose $\V$ so that $\ol{\V}$
is disjoint from $\sigma_1(\ol{\V})$.  Therefore in this case,
$A$ is a continuous function.
 
Choose $h \in \rC_0(Y)$ such that $h|_{\gamma_s(\ol{\V})} = 1$
and $\|h\|_\infty = 1$.
Define $\psi_j(z)$ to be the $1,2$ entry of $\Phi(x_0, z)( \ft_j h)$;
and set $\Psi(z) = (\psi_1(z),\dots,\psi_n(z))$.
Then $\Psi$ is an analytic function from $\bC^k$ into $\bC^n$.

\textit{We claim that $\psi_j(z)=0$ for $j>l$.}

Indeed, since $j>l$, the map
$\gamma^{-1} \tau_j \gamma$ is not in $\F$.
Hence there exists a net $(x_ \lambda)_{\lambda \in \Lambda}$
in $\ol{\V}$ converging to $x_0$ so that
$\gamma^{-1}_{s} \tau_j \gamma_{s}(x_\lambda) \ne  \sigma_1(x_\lambda)$
for all $\lambda \in \Lambda$.
By Lemma~\ref{L:xtoy}, $\Phi(x_\lambda,z) (\ft_j h)$ is diagonal
for all $\lambda$ in $\Lambda$.

First consider the case when $\sigma_1(x_0) \neq x_0$.
Then $A(x,z)$ is continuous, and so $\Phi(x,z)$ is point--norm continuous.
Taking limits, we conclude that $\Phi(x_0,z)(\ft_j h)$ is diagonal;
whence $\psi_j(z)=0$.

Now consider the case $\sigma_1(x_0) =  x_0$. Recall that in this case,
$\Phi(x_0,z)$ has range in $\A(E_{12})$; so that the diagonal part
consists of scalars. Fix $z \in \bC^k$.
Since
\[
 \max \{ \| A(x_\lambda, z) \| , \|A(x_\lambda, z)^{-1}\|\} \le 1 + \| \gamma_r \| \|z\| ,
\]
we may pass to a subnet if necessary so that
$\lim_{\Lambda} A(x_\lambda, z) = A(z)$ exists in $\fT_2^{-1}$.
Since $\Phi_0$ is point--norm continuous and $A(x_0,z) = I_2$,
\begin{align*}
 \lim_{\lambda \in \Lambda} \Phi(x_\lambda,z)(\ft_j h)
 &= \lim_{\lambda \in \Lambda} A(x_\lambda, z) \Phi_0(x_\lambda,z)(\ft_j h)A(x_\lambda, z)^{-1}
 \\&= A(z) \Phi_0(x_0,z)(\ft_j h) A(z)^{-1} .
\end{align*}
Therefore $A(z) \Phi_0(x_0,z)(\ft_j h) A(z)^{-1}$ is diagonal, and hence scalar.
So $\Phi(x_0,z)(t_jh)$ is scalar and $\psi_j(z) = 0$, which proves the claim.

The function $\Psi$ can now be considered as an analytic function
from $\bC^k$ into $\bC^l$.
Observe that $\Psi(0)=0$.
By Proposition \ref{zeroset} there exists $z_0\neq 0$ for which $\Psi(z_0)=0$.
Then $\Phi(z_0)$ is diagonal, and thus is not a nest representation.
This is a contradiction which proves the theorem.
\end{proof}

\section{The Converse}\label{S:converse}

In this section, we consider to what extent the converse is valid.
It turns out that for the tensor algebra, we are able to establish
this converse in many situations and for various notions of isomorphism.
We offer a plausible reason
for our conjecture that the converse holds in complete generality.
On the other hand, we will show by example that, unlike the tensor algebra situation,
the converse fails for
the semicrossed product if one considers completely isometric isomorphisms.

\begin{eg}\label{Ex:graph}
\textbf{The tensor algebra for a discrete set.}
Suppose that $X$ is a countable discrete set.
Then $\rC_0(X)$ is just $c_0(X)$.
In particular, the idempotents $p_x$ corresponding to characteristic
functions of each point in $X$ all belong to $\rC_0(X)$.
Therefore it is easy to see that $\A(X,\sigma)$ is generated by
$\{ p_x : x \in X\}$ and the partial isometries $\{ \fs_i p_x : x \in X,\, 1 \le i \le n\}$.
These partial isometries have domain $p_x$ and the complete set has
pairwise orthogonal ranges.
Evidently there is no way to determine any particular order on the
$n$ maps $\fs_1 p_x,\dots,\fs_n p_x$ independent of $x$.

There is a graph $G_\sigma$ associated to $(X,\sigma)$ with
vertices indexed by $X$, and $n$ edges with source $x$
and ranges $\sigma_1(x),\dots,\sigma_n(x)$.
This is the graph of $\sigma$ without the labels.
Evidently piecewise conjugacy in this context allows
arbitrary assignment of labels.
Therefore this graph represents a complete invariant up to
piecewise conjugacy of the system.

It is not difficult now to see that $\A(X,\sigma)$ is just the
tensor algebra of the graph.
The C*-envelope is the Cuntz--Krieger algebra of the graph
\cite{FMR, KK2}.
\end{eg}

\begin{eg} \label{Ex:pwc but not iso}
\textbf{The semicrossed product for a two point set.}
Consider two systems on $X = \{1,2\}$.
Let $\sigma_1$ be the identity map,
and let $\sigma_2$ interchange the two points.
Our first system is $(X,\sigma_1,\sigma_2)$.
Also let $\tau_1 = 1$ and $\tau_2 = 2$ be constant maps.
The second system is $(X,\tau_1,\tau_2)$.
We have
\[
 \sigma_1(1) = \tau_1(1),\ \sigma_2(1) = \tau_2(1) \qand
 \sigma_2(2) = \tau_1(2),\ \sigma_1(2) = \tau_2(2) .
\]
So these two systems are piecewise conjugate.
We will show that the semicrossed products of these
two systems are not isomorphic.

First consider $\rC(X) \times_\tau \Ftwo$.
This is generated by two complementary projections
$P_1$ and $P_2 = P_1^\perp$ and two isometries
$S_1$ and $S_2$ with ranges $P_1$ and $P_2$ respectively.
This is because in every full dilation, the range of $S_i$
will be the spectral projection for the range of $\tau_i$,
namely $\{i\}$.
This coincides with the tensor algebra $\A(X,\tau)$,
which by the preceding Example is just the tensor
algebra of the graph.
In particular, one can recognize within $\A(X,\tau)$
that there are proper isometries with complementary ranges.
It is easy to see that the C*-envelope is the Cuntz algebra $\O_2$.

Now consider $\rC(X) \times_\sigma \Ftwo$.
This algebra is generated by two complementary projections
and two unitaries:
\[
 P_1 = \begin{bmatrix}I&0\\0&0\end{bmatrix},\
 P_2 = \begin{bmatrix}0&0\\0&I\end{bmatrix},\
 U_1 = \begin{bmatrix}U_{11}&0\\0&U_{22}\end{bmatrix}
 \AND
U_2 = \begin{bmatrix}0&U_{12}\\U_{21}&0\end{bmatrix} .
\]
There are no relations on the $U_{ij}$ except for the fact
that there is no canonical identification between $P_1\H$
and $P_2\H$.  So we can take this identification to be
given by $U_{21}$.  That is, we have $U_{21}=I$ and
the other three unitaries are in general position.

The C*-algebra generated in this way is 
$\fM_2(\ca(\bF_3))$, the $2\times 2$ matrices
over the full group C*-algebra of the free group on 3 generators.
A representation of this C*-algebra is provided by arbitrary
choices for the unitaries $U_{11}$, $U_{22}$ and $U_{21}$.
It is clear that any such representation restricted to
$\rC(X) \times_\sigma \Ftwo$ is maximal, and thus 
factors through the C*-envelope.
Hence $\fM_2(\ca(\bF_3))$ is the C*-envelope.
This C*-algebra is finite, and in particular contains
no proper isometries. 
It follows that the two semicrossed products, $\rC(X) \times_\sigma \Ftwo$
and $\rC(X) \times_\tau \Ftwo$, are not completely isometrically isomorphic.
\end{eg}

We now restrict our attention
to the tensor algebra case.
We prove various partial converses.
These may be summarized as follows:

\begin{thm}\label{T:partial converse}
Suppose that at least one of the following holds:
\begin{itemize}
\item $n \le 3$, or
\item $X$ has covering dimension at most 1, or
\item $\Z(X,\sigma) = \{ x : |\sigma(x)| < n \}$ has no interior.
\end{itemize}
Then the following are equivalent:
\begin{enumerate}
\item $(X,\sigma)$ and $(Y,\tau)$ are piecewise topologically conjugate.
\item $\A(X,\sigma)$ and $\A(Y,\tau)$ are isomorphic.
\item $\A(X,\sigma)$ and $\A(Y,\tau)$ are completely isometrically isomorphic.
\end{enumerate}
\end{thm}

The evidence strongly suggests the following:

\begin{conj}\label{Conj:iso}
The tensor algebras $\A(X,\sigma)$ and $\A(Y,\tau)$ are isomorphic
if and only if the systems $(X,\sigma)$ and $(Y,\tau)$ are piecewise
topologically conjugate; and in this case, they are completely
isometrically isomorphic.
\end{conj}

We will now prove the various pieces of Theorem~\ref{T:partial converse}
and provide more evidence for our Conjecture~\ref{Conj:iso}.
 
When piecewise conjugacy reduces to conjugacy, we obtain the following
classification for both the tensor algebra and the semicrossed product.
This is the best result we can offer for the semi-crossed product.

\begin{cor}  \label{last}
Let $(X, \sigma)$ and $(Y, \tau)$ be multivariable dynamical systems,
and assume that $X$ is connected and  $\Z(X, \sigma)$ has empty interior.
Then $\A$ and $\B$ are isomorphic if and only if the systems
$(X, \sigma)$ and $(Y, \tau)$ are conjugate.
\end{cor}

In the totally disconnected case (dimension 0), the situation is straightforward.
The reader should compare this with \cite{Ion}.

\begin{cor}\label{C:total_disc}
Assume that $X$ is totally disconnected.
Then the tensor algebras $\A(X,\sigma)$ and $\A(Y,\tau)$
are isomorphic if and only if
$(X,\sigma)$ and $(Y,\tau)$ are piecewise topologically conjugate.
In this case, the algebras are completely isometrically isomorphic.
\end{cor}

\begin{proof}
One direction is provided by Theorem~\ref{class}.
For the converse, assume that $(X,\sigma)$ and $(Y,\tau)$
are piecewise topologically conjugate.
By Proposition~\ref{P:total_disc}, there
is a partition of $X$ into clopen sets $\{ \V_\alpha : \alpha\in S_n \}$
so that  for each $\alpha \in S_n$,
\[
 \gamma^{-1} \tau_i \gamma|_{\V_\alpha} = \sigma_{\alpha(i)} |_{\V_\alpha} .
\]

We may assume that $\A(X,\sigma)$ is contained in a universal operator algebra
generated by $\rC_0(X)$ and $n$ isometries $\fs_1,\dots,\fs_n$
satisfying the covariance relations.  Likewise, $\A(Y,\tau)$ is contained
in the corresponding algebra generated by $\rC_0(Y)$ and
$\ft_1,\dots,\ft_n$.
Define a covariant representation of $(Y,\tau)$ by
\begin{align*}
 \phi(f) &= f \circ \gamma \\
 \phi(\ft_i) &= \sum_{\alpha \in S_n} \fs_{\alpha(i)} \upchi_{\V_\alpha}
 \qfor 1 \le i \le n.
\end{align*}

To see that $\big[ \phi(\ft_i)\ \dots \  \phi(\ft_n) \big]$ is a row isometry,
compute
\begin{align*}
  \phi(\ft_j)^*  \phi(\ft_i) &= \sum_{\alpha,\beta \in S_n}
   \upchi_{\V_\beta} \fs_{\beta(j)}^* \fs_{\alpha(i)} \upchi_{\V_\alpha} \\
 &=  \sum_{\alpha,\beta \in S_n} \delta_{\beta(j),\alpha(i)}  \upchi_{\V_\beta} \upchi_{\V_\alpha} \\
 &= \sum_{\alpha \in S_n} \delta_{\alpha(j),\alpha(i)} \upchi_{\V_\alpha} = \delta_{j,i} I .
\end{align*}
So $\phi(t_i)$ are isometries with pairwise orthogonal ranges.

Next observe that $\phi$ satisfies the covariance relations.
\begin{align*}
 \phi(f) \phi(t_i) &= (f \circ\gamma) \sum_{\alpha \in S_n} \fs_{\alpha(i)} \upchi_{\V_\alpha} \\
 &= \sum_{\alpha \in S_n} \fs_{\alpha(i)} \upchi_{\V_\alpha}  (f\circ\gamma \sigma_{\alpha(i)}) \\
 \intertext{and since $\gamma\sigma_{\alpha(i)}=\tau_i \gamma$ on $\V_\alpha$,}
 &= \sum_{\alpha \in S_n} \fs_{\alpha(i)} \upchi_{\V_\alpha} (f \circ \tau_i \gamma)
  = \phi(\ft_i)  \phi(f\circ \tau_i) .
\end{align*}

Therefore $\phi$ extends to a completely contractive representation
of $\A(Y,\tau)$ given by
\[ \phi \big( \sum_{w \in \Fn} \ft_w f_w \big) = \sum_{w \in \Fn} \phi(\ft_w) (f_w \circ \gamma) \]
where $w=i_k\dots i_1$ and $\phi(\ft_w) = \phi(\ft_{i_k}) \dots \phi(\ft_{i_1})$.
It is evident that this maps $\A(Y,\tau)$ into $\A(X,\sigma)$.

However the relations between $\sigma$ and $\tau$ can be reversed to obtain
a completely contractive map from $\A(X,\sigma)$ into $\A(Y,\tau)$.
A little thought shows that this map is the inverse of $\phi$,
verifying that $\phi$ is a completely isometric isomorphism.
\end{proof}

Let us define the graph of $\sigma$ to be
\[ G(\sigma) = \{ (x,\sigma_i(x)) : x \in X,\ 1 \le i \le n \} \]
considered as a subset of $X \times X$.

\begin{cor}\label{C:disjoint}
If the maps $\sigma_i$ for $1 \le i \le n$ have disjoint graphs,
then $\A(X,\sigma)$ is isomorphic to $\A(Y,\tau)$ if and
only if there is a homeomorphism $\gamma$ of $X$ onto $Y$
which implements a homeomorphism between the graphs
of $\sigma$ and $\tau$.
In this case, the algebras are completely isometrically isomorphic.
\end{cor}

\begin{proof}
By Theorem~\ref{class}, an isomorphism between $\A(X,\sigma)$
and $\A(Y,\tau)$ yields a piecewise topological conjugacy
via a homeomorphism $\gamma$.
In particular, this implements a homeomorphism from the
graph $G(\sigma)$ onto $G(\tau)$ via the map $\gamma \times \gamma$.

Conversely suppose that $\gamma$ is a homeomorphism of $X$ onto $Y$
so that $\gamma \times \gamma$ carries $G(\sigma)$ onto $G(\tau)$.
For simplicity of notation, we may suppose that $Y=X$ and $\gamma = \id$.
The fact that $\sigma_i$ have disjoint graphs means that
$|\sigma(x)| = n$ for each $x \in X$.
Since $\tau(x) = \sigma(x)$, there is a unique permutation $\alpha_x \in S_n$
so that $\tau_i(x) = \sigma_{\alpha_x(i)}(x)$ for $1 \le i \le n$.
A simple argument using the continuity of $\sigma$ and $\tau$
shows that the map $\alpha$ taking $x \in X$ to $\alpha_x \in S_n$
is continuous.
Therefore the sets $\V_\alpha = \{ x : \alpha_x = \alpha \}$
yields a partition of $X$ into clopen sets with the property that
\[
 \tau_i|_{\V_\alpha} = \sigma_{\alpha(i)}|_{\V_\alpha}
 \qfor \alpha \in S_n \AND 1 \le i \le n .
\]
The proof is now completed as for Corollary~\ref{C:total_disc}.
\end{proof}

When $n=2$, piecewise conjugacy can be described in simple terms.
This also leads to a fairly straightforward converse which points
the way to the issues in higher dimensions.

\begin{thm}\label{T:disjoint_closures}
Let $n=2$ and suppose that $(X,\sigma_1,\sigma_2)$
and $(Y,\tau_1,\tau_2)$ are two dynamical systems.
The following are equivalent:
\begin{enumerate}
\item $(X,\sigma)$ and $(Y,\tau)$ are piecewise topologically conjugate.
\item $\A(X,\sigma)$ and $\A(Y,\tau)$ are isomorphic.
\item $\A(X,\sigma)$ and $\A(Y,\tau)$ are completely isometrically isomorphic.
\item There is a homeomorphism $\gamma$ of $X$ onto $Y$ so that
\begin{enumerate}
\item $\{ \gamma \sigma_1(x),  \gamma \sigma_2(x) \} =
 \{ \tau_1 \gamma (x),  \tau_2 \gamma (x) \}$ for each $x \in X$, and
\item $X_i = \{ x \in X : \gamma \sigma_1(x) \ne \tau_i \gamma (x) \}$, $i=1,2$, have disjoint closures.
\end{enumerate}
\end{enumerate}
\end{thm}

\begin{proof}
Clearly (3) implies (2); and (2) implies (1) by Theorem~\ref{class}.
Suppose that (1) holds for a homeomorphism $\gamma$ and an
open cover  $\{ \U_\id, \U_{(12)} \}$.
That is, $\tau_i \gamma|_{ \U_\id} = \gamma \sigma_i|_{ \U_\id}$
and $\tau_i \gamma|_{ \U_{(12)} } = \gamma \sigma_{i'}|_{ \U_{(12)} }$
for $i=1,2$ and $\{ i,i'\} = \{1,2\}$.
In particular, (4a) holds for every $x \in X$.
Moreover, $\ol{X_1} \subset \U_\id^c$ and $\ol{X_2} \subset \U_{(12)}^c$.
Consequently,
\[ \ol{X_1} \cap \ol{X_2} \subset \big( \U_\id \cup \U_{(12)} \big)^c = \mt .\]
This establishes (4b).

So we now assume that (4) holds.  To simplify the notation,
we may assume that $Y=X$ and $\gamma = \id$.
Let $(\pi,S_1,S_2)$ be a faithful representation of the
covariance relations for $(X,\sigma)$.
Since $\ol{X_1}$ and $\ol{X_2}$ are disjoint,
there is a continuous function $h \in \rC_b(X)$
such that $h|_{\ol{X_1}} = 0$, $h|_{\ol{X_2}}=\frac\pi 2$
and $h$ takes real values in $[0,\frac\pi 2]$ everywhere.
Let $\ol{\pi}$ denote the extension of $\pi$ to the
bounded Borel functions on $X$.
Define
\[
 T_1 = S_1 \ol{\pi}(\sin h) + S_2  \ol{\pi}(\cos h) \qand
 T_2 = S_1  \ol{\pi}(\cos h) - S_2  \ol{\pi}(\sin h) .
\]

First observe that $T_i$ are isometries with orthogonal range.
For example,
\[ T_i^*T_1 =  \ol{\pi}(\sin h) S_1^*S_1  \ol{\pi}(\sin h) +
 \ol{\pi}(\cos h) S_2^*S_2 \ol{\pi}(\cos h) = I \]
and
\[ T_1^*T_2 = \ol{\pi}(\sin h) S_1^*S_1 \ol{\pi}(\cos h) -
 \ol{\pi}(\cos h) S_2^*S_2 \ol{\pi}(\sin h) = 0 .\]
Next verify the covariance relations $\pi(f)T_i = T_i \pi(f\circ \tau_i)$:
\begin{align*}
 \pi(f) T_i &= \pi(f) S_1 \ol{\pi}(\sin h) + \pi(f) S_2  \ol{\pi}(\cos h)\\
 &= S_1 \ol{\pi}(\sin h) \pi(f\circ\sigma_1) + S_2  \ol{\pi}(\cos h) \pi(f\circ\sigma_2)\\
 \intertext{and since $\sigma_1=\tau_1$ on $h^{-1}((0,\frac\pi 2])$
 and $\sigma_2=\tau_1$ on $h^{-1}([0,\frac\pi 2))$,}
 &= \big( S_1 \ol{\pi}(\sin h) + S_2  \ol{\pi}(\cos h) \big) \pi(f\circ\tau_1) \\
 &= T_i \pi(f\circ \tau_i) .
\end{align*}

Next observe that $\A(X,\sigma)$ is generated by $\rC_0(X)$ and
$T_i\rC_0(X)$ for $i=1,2$.  This is because
\[
 S_1 = T_1 \ol{\pi}(\sin h) + T_2  \ol{\pi}(\cos h) \qand
 S_2 = T_1  \ol{\pi}(\cos h) - T_2  \ol{\pi}(\sin h) .
\]
Multiplying on the right by $\pi(f)$ for any $f \in \rC_0(X)$
yields the corresponding fact for $\A(X,\sigma)$.

Now exactly the same procedure works in a faithful representation
of $\A(X,\tau)$.
Let $\ft_1$ and $\ft_2$ denote the generators of $\A(X,\tau)$.
Since $\A(X,\sigma)$ contains a representation of the covariance
relations for $(X,\tau)$, there is a homomorphism of
$\A(X,\tau)$ onto $\A(X,\sigma)$ that
takes $f$ to $\pi(f)$ and $\ft_i f$ to $T_i f$.
It is then clear that $a_1 := \ft_1 \ol{\pi}(\sin h) + \ft_2  \ol{\pi}(\cos h)$
is sent to $S_1$ and $a_2 := \ft_1  \ol{\pi}(\cos h) - \ft_2  \ol{\pi}(\sin h)$
is sent to $S_2$.
Likewise there is an algebra homomorphism in the reverse direction
which is evidently the inverse map on the generators.
Therefore these algebras are isomorphic.
Moreover the maps in both directions are complete contractions,
and thus are completely isometric.
\end{proof}

The analysis leads to a technical conjecture which we can verify in
low dimensions that would suffice to solve our Conjecture~\ref{Conj:iso}.

\begin{conj}\label{Conj:simplex}
Let $\Pi_n$ be the $n!$-simplex with vertices indexed by $S_n$.
Then there should be a continuous function $u$ of $\Pi_n$  into $U(n)$
so that:
\begin{enumerate}
\item  each vertex is taken to the corresponding permutation matrix,
\item for every pair of partitions $(A,B)$ of the form
\[
 \{1,\dots,n\} = A_1 \dot\cup \dots \dot\cup A_m = B_1 \dot\cup \dots \dot\cup B_m,
\]
where $|A_s|=|B_s|,\ 1\le s \le m$,
let
\[ \P(A,B) = \{\alpha \in S_n : \alpha(A_s)=B_s, 1\le s \le m \} .\]
If $x = \sum_{\alpha \in \P(A,B)} x_\alpha \alpha$, then
the non-zero matrix coefficients of $u_{ij}(x)$ are
supported on $\bigcup_{s=1}^m B_s \times A_s$.
We call this the \textit{block decomposition condition}.
\end{enumerate}
\end{conj}

We establish this conjecture for $n=2$ and $3$.
We will then demonstrate how this conjecture solves
our problem, and in particular provides the solution when $n=3$.

\begin{prop}\label{P:conj23}
Conjecture~$\ref{Conj:simplex}$ is valid for $n=2,3$.
Moreover there is a function on the $1$-skeleton of $\Pi_n$
satisfying the conditions of this conjecture for every $n$.
\end{prop}

\begin{proof} For each permutation $\alpha \in S_n$, one can choose a Hermitian
matrix $A_\alpha$ so that $U_\alpha = \exp(iA_\alpha)$ as follows:
Decompose $\alpha$ into cycles, and select a logarithm for
each cycle with arguments in $[-\pi,\pi]$.  The eigenvalues come in
conjugate pairs except for $\pm1$.  Make the choice
between $\pm\pi$ for the eigenvalue $-1$ alternately; so that
$\|A_\alpha\| \le \pi$ and $\Tr A_\alpha \in \{0,\pi\}$.
Choosing $A_\alpha$ in this way ensures that it respects the
block diagonal structure of $U_\alpha$ coming from the cycle decomposition.

To describe such a function on the $1$-skeleton of $\Pi_n$,
just order the elements of $S_n$.
For future use, we insist that the even permutations precede the odd
permutations in this order.
Then if $\alpha < \beta$ in this order, define
\[ u((1-t)\alpha + t\beta) = U_\alpha \exp(itA_{\alpha^{-1}\beta}) .\]
This works because if $\alpha,\beta \in P(A,B)$, then $\alpha^{-1}\beta \in P(A,A)$.
Hence $A_{\alpha^{-1}\beta}$ respects the block structure $A$,
so that $\exp(itA_{\alpha^{-1}\beta})$ does also.
Consequently $U_\alpha \exp(itA_{\alpha^{-1}\beta})$ respects the $A,B$ block decomposition.
This $1$-skeleton argument contains the $n=2$ case.

Observe that for any $n$, the function
$u(\sum t_\alpha \alpha) =\exp(i \sum t_\alpha A_\alpha)$
maps $\Pi_n$ into $U(n)$ and $u(\alpha) = U_\alpha$.
n general, it does not satisfy the block decomposition condition.

Now consider $n=3$. Any non-trivial block decomposition of $3\times 3$
is given by a vertex $(i,j)$; namely
\[\{ A_1=\{i\},\, A_2 = A_1^c,\, B_1=\{j\},\, B_2=B_1^c \} .\]
There are two permutations respecting this block decomposition,
say $\alpha_{ij}$ and $\beta_{ij}$ where $\alpha_{ij}$ is even and $\beta_{ij}$ is odd.
To find the function $u$ satisfying Conjecture~\ref{Conj:simplex},
we start with the function constructed in the previous paragraph
and modify it on the nine edges to obtain the block decomposition condition.

Take the $6$-simplex $\Pi_3$ and glue a half-disk $D_{ij}$ to each of the
nine edges$[\alpha_{ij}, \beta_{ij}]$.
On the semicircular edge $\Gamma_{ij}[0,1]$ between $\alpha_{ij}$ and $\beta_{ij}$,
define $u$ as in the second paragraph:
\[ u(\Gamma_{ij}(t)) = U_{\alpha_{ij}} \exp(itA_{\alpha_{ij}^{-1}\beta_{ij}}) .\]
Recall that on the straight edge $[\alpha_{ij},\beta_{ij}]$, $u$ is define as
\[ u((1-t)\alpha_{ij} + t\beta_{ij}) = \exp(i(1-t)A_{\alpha_{ij}} + itA_{\beta_{ij}}) .\]
The determinant along these two paths is respectively
\[ \det(U_{\alpha_{ij}}) \exp(it \Tr A_{\alpha_{ij}^{-1}\beta_{ij}}) = e_{it\pi} \]
and
\[ \exp(i(1-t) \Tr A_{\alpha_{ij}} + it \Tr A_{\beta_{ij}}) = e_{it\pi} \]
because $\Tr A_{ij}=0$ and
$\Tr A_{\beta_{ij}} = \Tr A_{\alpha_{ij}^{-1}\beta_{ij}} = \pi$.

The issue is now to extend the definition of $u$ to each half-disk $D_{ij}$
in a continuous way.
The key fact that we need is that $SU(3)$ is simply
connected \cite[Theorem~II.4.12]{MiTo}.
(Indeed $SU(n)$ is simply connected for all $n\ge1$.)
The determinant is an obstruction to $U(n)$ being simply connected.
But by ensuring that our functions have domain in the set of unitaries
with determinant in the upper half plane means that we remain in a
simply connected set.
Hence there is a homotopy between the edge $u([\alpha_{ij},\beta_{ij}])$
and $u(\Gamma_{ij})$ which enables the extension of the definition of $u$.
Finally observe that there is a continuous function $h$ of $\Pi_3$
onto $\Pi_3 \bigcup_{i,j=1}^3 D_{ij}$ which takes the edges
$[\alpha_{ij},\beta_{ij}]$ onto $\Gamma_{ij}$.
The composition $v = u\circ h$ is the desired function.
\end{proof}

We can now show how to modify the proof of Theorem~\ref{T:disjoint_closures}
to work for any $n$ for which we have Conjecture~\ref{Conj:simplex}.
Recall that the covering dimension of a topological space $X$ is the
smallest integer $k$ so that every open cover can be refined so that
each point is covered by at most $k+1$ points.

\begin{thm}\label{T:main_converse}
If Conjecture~$\ref{Conj:simplex}$ is correct for some value of $n$, then
for two paracompact dynamical systems $(X,\sigma)$ and $(Y,\tau)$, where
$\sigma = \{\sigma_1,\dots,\sigma_n\}$ and $\tau=\{\tau_1,\dots,\tau_n\}$,
the following are equivalent:
\begin{enumerate}
\item $(X,\sigma)$ and $(Y,\tau)$ are piecewise topologically conjugate.
\item $\A(X,\sigma)$ and $\A(Y,\tau)$ are isomorphic.
\item $\A(X,\sigma)$ and $\A(Y,\tau)$ are completely isometrically isomorphic.
\end{enumerate}
If Conjecture~$\ref{Conj:simplex}$ is valid on the $k$-skeleton of $\Pi_n$,
then the theorem still holds if the covering dimension of $X$ is at most $k$.
\end{thm}

\begin{proof}
Clearly (3) implies (2); and (2) implies (1) by Theorem~\ref{class}.
Suppose that (1) holds for a homeomorphism $\gamma$ and an
open cover $\{ \U_\alpha : \alpha \in S_n \}$.
To simplify notation, we may identify $Y$ with $X$ via $\gamma$,
so that we have $Y=X$ and $\gamma = \id$.
If the covering dimension of $X$ is $k$, then the cover
$\{ \U_\alpha : \alpha \in S_n \}$ can be refined so that
each point of $X$ is contained in at most $k+1$ open sets.
Let $\{ g_\alpha : \alpha \in S_n \}$ be a partition of unity relative to this cover.
This induces a map $g$ of $X$ into the simplex $\Pi_n$ given by
$g(x) =\big( g_\alpha(x) \big)_{\alpha\in S_n}$.
In the case of covering dimension $k$, the image lies in the $k$-skeleton.

Our hypothesis is that there is a continuous function $u$ from $\Pi_n$
or its $k$-skeleton into the unitary group $U(n)$ satisfying the block
decomposition condition.  Let $v = u \circ g$ map $X$ into $U(n)$.
Now if $x\in X$, there is a minimal partition
\[
 \{1,\dots,n\} = A_1 \dot\cup \dots \dot\cup A_m
 = B_1 \dot\cup \dots \dot\cup B_m
\]
into maximal subsets and an open neighbourhood $\U$ of $x$
so that $\sigma_i |_\U = \tau_j |_\U$ for $i\in A_s$ and $j\in B_s$,
$1 \le s \le m$.
The permutations $\alpha$ for which $g_\alpha(x) \ne 0$
respect this block structure.
Hence so does the map $v$.
This will ensure that in our construction below, we will always intertwine
functions that agree on a neighbourhood of each point.

Let $v_{ij}$ be the matrix coefficients of $v$.
Define operators in $\A(X,\sigma)$ by $T_i = \sum_{j=1}^n \fs_j v_{ij}$.
Then since the $\fs_j$'s have pairwise orthogonal range,
\[
 T_k^* T_i = \sum_{j=1}^n \ol{v_{kj}} v_{ij} = \delta_{ki} I .
\]
Hence the $T_i$'s are isometries with pairwise orthogonal ranges.
To check the covariance relations, observe that if $v_{ij}(x) \ne0$,
then $\tau_i$ and $\sigma_j$ agree on a neighbourhood of $x$.
That is, $v_{ij}  (f\circ\sigma_j) = v_{ij}  (f\circ\tau_i)$ for all $i,j$.
Therefore
\begin{align*}
 f T_i &= f \sum_{j=1}^n \fs_j v_{ij}
 = \sum_{j=1}^n \fs_j v_{ij}  (f\circ\sigma_j)\\
 &= \sum_{j=1}^n \fs_j v_{ij}  (f\circ\tau_i)
 = T_i (f\circ \tau_i) .
\end{align*}

Next observe that $\A(X,\sigma)$ is generated by $C_0(X)$
and $T_i C_0(X)$ for $1 \le i \le n$.  This is because for $1\le k \le n$,
\begin{align*}
 \sum_{i=1}^n T_i\ol{v_{ik}} f
 &= \sum_{i=1}^n \sum_{j=1}^n \fs_j v_{ij} \ol{v_{ik}} f  
 = \sum_{j=1}^n \fs_j f \sum_{i=1}^n v_{ij} \ol{v_{ik}} = \fs_k f .
\end{align*}
Therefore there is a completely contractive homomorphism of
$\A(Y,\tau)$ onto $\A(X,\sigma)$ sending $\ft_i$ to $T_i$ for
$1 \le i \le n$ which is the identity on $C_0(Y)=C_0(X)$.
Likewise there is a completely contractive homomorphism of
$\A(X,\sigma)$ onto $\A(Y,\tau)$ which is the inverse on the
generators $\fs_j$.
Consequently these maps are completely isometric isomorphisms.
\end{proof}

As Proposition~\ref{P:conj23} establishes
Conjecture~\ref{Conj:simplex} for $n=3$, we obtain
the immediate consequence:

\begin{cor}\label{C:n=3}
Suppose that $n=3$.
Then for two dynamical systems $(X,\sigma)$
and $(Y,\tau)$, the following are equivalent:
\begin{enumerate}
\item $(X,\sigma)$ and $(Y,\tau)$ are piecewise topologically conjugate.
\item $\A(X,\sigma)$ and $\A(Y,\tau)$ are isomorphic.
\item $\A(X,\sigma)$ and $\A(Y,\tau)$ are completely isometrically isomorphic.
\end{enumerate}
\end{cor}

The dimension $0$ case corresponds to totally disconnected spaces.
So this result subsumes Corollary~\ref{C:total_disc}.
As compact subsets of $\bR$ have covering dimension $1$,
we obtain:

\begin{cor}\label{C:1dim}
Suppose that $X$ has covering dimension $0$ or $1$.
In particular, this holds when $X$ is a compact subset of $\bR$.
Then for two dynamical systems $(X,\sigma)$
and $(Y,\tau)$, the following are equivalent:
\begin{enumerate}
\item $(X,\sigma)$ and $(Y,\tau)$ are piecewise topologically conjugate.
\item $\A(X,\sigma)$ and $\A(Y,\tau)$ are isomorphic.
\item $\A(X,\sigma)$ and $\A(Y,\tau)$ are completely isometrically isomorphic.
\end{enumerate}
\end{cor}

\chapter{Semisimplicity} \label{C:semisimple}
\section{Wandering sets and Recursion}\label{S:wander}

In this section, we examine some topological issues which are needed
in the next section on semisimplicity.

\begin{defn}
Let $(X , \sigma)$ be a multivariable dynamical system.
An open set $U \subseteq X$ is said to be \textit{$(u,v)$-wandering} if
\[
\sigma_{u w v}^{-1}(U) \cap U = \mt \qforal w \in \Fn .
\]
If $U$ is a $(u,v)$-wandering set for some $u,v$ in $\Fn$,
we say that $U$ is a \textit{generalized wandering set}.
We may write \textit{$v$-wandering} instead of $(\mt,v)$-wandering.
\end{defn}

Observe that if some $\sigma_i$ is not surjective, then the
set $U = X\bsl \sigma_i(X)$ is $(i,\mt)$-wandering.
If one is interested in when there are no wandering sets,
one should restrict to the case in which each $\sigma_i$ is surjective.
In this case, whenever $U$ is $(u,v)$-wandering,
the set $\sigma_{u}^{-1}(U)$ is $(\mt,vu)$-wandering.
This open set is non-empty because every $\sigma_u$ is surjective.
In general, the additional generality seems to be needed.

\begin{defn}
Let $(X , \sigma)$ be a multivariable dynamical system.
Given $u, v \in \Fn$, we say that a point $x \in X$ is
\textit{$(u,v)$-recurrent} if for every neighbourhood $U \ni x$,
there is a $w \in \Fn$ so that $\sigma_{uwv}(x) \in U$.
We will write \textit{$v$-recurrent} instead of $(\mt,v)$-recurrent.
\end{defn}

Again, if each $\sigma_i$ is surjective, we see that
when $x$ is $(u,v)$-recurrent, then any point
$y$ in the non-empty set $\sigma_{u}^{-1}(x)$ is $(\mt,vu)$-recurrent.

The following proposition explains the relationship
between recursion and wandering sets in metric spaces.

\begin{prop}\label{P:wander}
Let $(X,\sigma)$ be a metrizable multivariable dynamical system.
There is no $(u,v)$-wandering set $U \subset X$ if and only if
the set of $(u,v)$-recurrent points is dense in $X$.
\end{prop}

\begin{proof}
If the $(u,v)$-recurrent points are dense in $X$, then any
open set $U \subset X$ contains such a point, say $x_0$.
Thus there is a word $w \in \Fn$ so that
$\sigma_{uwv}(x_0) \in U$.
Therefore $\sigma_{uwv}^{-1}(U) \cap U \ne \mt$.
So there are no $(u,v)$-wandering sets.

Conversely, suppose that there are no $(u,v)$-wandering sets.
Let $\rho$ be the metric on $X$.
For each $x \in X$, define
\[
 \delta(x) = \delta_{u,v}(x)
 := \inf_{w \in \Fn} \rho(x,\sigma_{uwv}(x)) .
 \]
Observe that $\delta$ is upper semicontinuous, i.e.\
$\{x \in X : \delta(x) < r \}$ is open for all $r \in \bR$.
Indeed, if $\delta(x_0) <r$ and $\rho(x_0,uwv(\sigma)(x_0)) = r - \ep$,
the continuity of our system means that for some $0 < \delta < \ep/2$,
$\rho(x,x_0) < \delta$ implies that
$\rho(\sigma_{uwv}(x),\sigma_{uwv}(x_0)) < \ep/2$.
So by the triangle inequality, $\rho(x,\sigma_{uwv}(x)) < r$;
whence $\delta(x) < r$.

Suppose that $U$ is a non-empty open set containing no
$(u,v)$-recurrent points.
Then $\delta(x) > 0$ for every $x \in U$.  So
\[ U = \bigcup_{n\ge1} \{ x \in U : \delta(x) \ge \frac1n \} .\]
This is a union of closed sets.  By the Baire Category Theorem,
there is an integer $n_0$ so that $\{ x \in U : \delta(x) \ge \frac1{n_0} \}$
has non-empty interior, say $V$.

Select a ball $U_0$ contained in $V$ of diameter less than $1/n_0$.
We claim that $U_0$ is $(u,v)$-wandering.
Indeed, for any $x \in U_0$ and any $w \in \Fn$,
$\rho(x,\sigma_{uwv}(x)) \ge 1/n_0$ and thus $uwv(\sigma)(x)$
lies outside of $U_0$.  That is,
$U_0 \cap \sigma_{uwv}^{-1}(U_0)  = \mt$.
\end{proof}

\begin{cor}\label{C:metric_wander}
In a metrizable multivariable dynamical system $(X,\sigma)$,
the following are equivalent:
\begin{enumerate}
\item there is no non-empty generalized wandering set in $X$

\item the set of $(u,v)$-recurrent points are dense for every $u,v \in \Fn$

\item each $\sigma_i$ is surjective and the
set of $v$-recurrent points is dense in $X$ for every $v \in \Fn$.
for every $v \in \Fn$.
\end{enumerate}
\end{cor}

\begin{proof} The equivalence of (1) and (2) is immediate
from the Theorem.
It was observed that (1) implies that each $\sigma_i$ is surjective;
and clearly (2) implies as a special case that the
set of $v$-recurrent points is dense in $X$ for every $v \in \Fn$.

Conversely, suppose that (3) holds, and fix $(u,v)$ and a point $x_0 \in X$.
By surjectivity, there is a point $y_0\in X$ so that $\sigma_u(y_0)=x_0$.
Given a neighbourhood $U_0 \ni x_0$, let $V_0 = \sigma_u^{-1}(U_0)$.
Select a $vu$-recurrent point $y \in V_0$; and let $x = \sigma_u(y)$.
If $U \ni x$ is any open set, let $V = \sigma_u^{-1}(U)$.
There is a word $w$ so that $\sigma_{wvu}(y) \in V$.
Therefore $\sigma_{uwv}(x) \in U$.
Thus $x$ is $(u,v)$-recurrent, and the set of such points
is dense in $X$. So (2) holds.
\end{proof}

\section{Semisimplicity}\label{S:semisimple}
 
We now turn to an analysis of the tensor algebra and the semicrossed product
to decide when they are semisimple.

An ideal $\J$ of an algebra $\A$ is said to be primitive if it is the kernel of an
algebraically irreducible representation. The intersection of all primitive ideals of $\A$ is
the Jacobson radical of $\A$. If the Jacobson radical of $\A$
happens to equal $\{ 0\}$, then $\A$ is said to be
\textit{semisimple}. If $\A$ is a Banach algebra, then the Jacobson
radical is closed since every primitive
ideal is the kernel
of some continuous representation of $\A$ on a Banach space. In fact an
element $ A \in \A$ belongs to the radical iff the spectral radious
of $BA$ vanishes for all $B \in \A$. This is a key property of the Jacobson radical
that is used in the sequel. Details can be found in \cite{Pal}.

We now apply the results of the previous section to characterize the semisimplicity
in our setting. In spite of the differences betwen the two
universal algebras, the answer turns out to be the same.

\begin{thm} \label{semisimple}
Let $(X,\sigma)$ be a multivariable dynamical system.
The following are equivalent:
\begin{enumerate}
\item $\A( X , \sigma)$  is semisimple.
\item $\rC_0(X) \times_\sigma \Fn$  is semisimple.
\item There are no non-empty generalized wandering sets.
\end{enumerate}
When $X$ is metrizable, these are also equivalent to
\begin{enumerate}
\setcounter{enumi}{3}
\item Each $\sigma_i$ is surjective and the $v$-recurrent points
are dense in $X$ for every $v\in \Fn$.
\end{enumerate}
\end{thm}

\begin{proof}
Assume first that there exists a nonempty open set $U \subset X$ and
$u, v \in \Fn$ so that
\[
\sigma_{uwv}^{-1}(U)  \cap U = \mt \qforal w \in \Fn.
\]
Let $h \ne 0$ be continuous function with support contained in $U$.
Then $(h\circ \sigma_{uwv}) h = 0$.

We will show that $\A(X , \sigma)$ is not semisimple.
Indeed, the (non-zero) operator $N= \fs_v h \fs_u$ generates a nilpotent ideal.
To see this, let $w \in \Fn$ and $f \in \rC_0(X)$.
Compute
\begin{align*}
N (\fs_w f ) N &=
\fs_v h \fs_{uw} f \fs_v h \fs_u     \\
&=\fs_{vuwv} (h\circ \sigma_{uwv}) h  (f \circ \sigma_v) \fs_u
 = 0 .
\end{align*}
Therefore $NAN=0$ for all $A \in \A(X,\sigma)$.
Hence the 2-sided ideal $\lip N \rip$ generated by $N$ is nilpotent of order 2;
and thus $\A(X,\sigma)$ is not semisimple.
The same calculation holds in the semicrossed product.
So both (2) and (3) imply (1).

Conversely, assume that there are no non-empty wandering sets.
We will show that both (2) and (3) hold.
As before, we will write $\A$ to denote either the tensor algebra or the
semicrossed product.
The Jacobson radical  of any Banach algebra is invariant under automorphisms.
In particular, both of our algebras have the automorphisms $\alpha_\lambda$
which send the generators $\fs_i$ to $\lambda_i \fs_i$ for each $\lambda\in\bT^n$.
Integration yields the expectations
\[ \Psi_\bk(a) = \int_{\bT^n} \alpha_\lambda(a) \ol{\lambda}^\bk \,d\lambda \]
for each $\ba \in \bN_0^n$ onto the polynomials spanned by
$\fs_w f$, where $w(\lambda) = \lambda^\bk$.
Consequently, these expectations map the radical into itself.

By way of contradiction, assume that $\rad \A$ contains non-trivial elements.
By the previous paragraph, $\rad \A$ will contain a non-zero element of the form
\[
 Y = \Phi_\ba(Y)
 = \sum_{w(\lambda) = \lambda^\bk} \fs_w h_w
 =: \sum_{j=1}^p \fs_{w_j} h_j
\]
for some $\ba \in \bN_0^n$.
Since $Y \ne 0$, we may suppose that $h_1 \ne 0$.
By multiplying $Y$ on the right by a function, we may suppose that
$h_1 \ge 0$ and $\|h_1\|>1$.

We will look for an element of the form
\[
 Q = \big( I + \sum_{k\ge1} 2^{-k}\fs_{v_k} \big) Y
    = \sum_{k\ge0} 2^{-k} \sum_{j=1}^p \fs_{w_{k,j}} h_j
\]
where $w_{k,j} = v_k w_j$.
The goal will be to show that $Q$ is not quasinilpotent, so contradicting
the fact that $Y$ is in the radical.

When $n=1$, $Y$ has only one term and the coefficients in $Q$
and all of its powers are non-negative functions.
Here there is no cancellation of terms, so the size of
a single term in the product $Q^m$ yields a useful lower bound
for the norm.

However when $n \ge 2$, this is not the case.
Instead, we can arrange to choose the words $v_k$
so that the distinct terms in a product $Q^m$ will
be distinct words in $\Fn$.
This will be accomplished if the word $v_k$ begins with the
sequence $t_k = 2 1^k 2$, the sequence consisting of
$k$ 1s separating two 2s.
For such $v_k$, consider the terms in a product $Q^m$.
They will have the form $\fs_u g_u$ where
$u$ is a product of $m$ words $u = w_{k_1,j_1} \dots w_{k_m,j_m}$.
The first few letters of $u$ will be $t_{k_1}$, uniquely
identifying $k_1$ and $w_{k_1,j_1}$ will be the first $|v_{k_1}|+|a|$
terms of $u$.  Peel that term off and repeat.
Since the product determines the terms, in order, it follows
that distinct products yield distinct words.
So again, there is no cancellation of terms.
Therefore the size of a single term in the product will yield
a lower bound for the norm of $Q$.

The plan of attack is to show that
\[ Q^{2^k-1} = 2^{-n_k}\fs_{u_k} g_k + \text{other terms} \]
where $\|g_k\| > 1$ and $n_k = 2^k-k-1$.
In computing the next power, we see that
\begin{align*}
 Q^{2^{k+1}-1} &= Q^{2^k-1} Q Q^{2^k-1} \\
 &= (2^{-n_k} \fs_{u_k} g_k)( 2^{-k} \fs_{w_{k,1}} h_1)( 2^{-n_k} \fs_{u_k} g_k)
      + \text{other terms}\\
 &= 2^{-2n_k-k} \fs_{u_k w_k u_k} (g_k\circ \sigma_{w_k u_k}) (h_1\circ \sigma_{u_k}) g_k
   + \text{other terms}\\
 &= 2^{-n_{k+1}} \fs_{u_{k+1}} g_{k+1} + \text{other terms} .
\end{align*}
Provided that $\|g_k\|>1$ for all $k \ge 1$, we obtain that
\[ \spr(Q) = \lim_{k\to\infty} \|Q^{2^k-1}\|^{1/(2^k-1)}
 \ge \lim_{k\to\infty} 2^{-\frac{2^k-k-1}{2^k-1}} = \frac12 > 0.
\]

The preceding calculation indicates where we should look.
One has formulae
\[ w_{k,1} = v_k w_1 ,\quad v_0=\mt, \qand u_{k+1} = u_k w_{k,1} u_k \qfor k\ge0 , \]
and $v_k$ begins with $t_k$ for $k\ge1$.
We choose $u_0=\mt$ so that $u_1 = w_1$ to start the induction.
We also have functions $g_1 = h_1$ and
\[ g_{k+1} = (g_k\circ \sigma_{w_{k,1} u_k}) (h_1\circ \sigma_{u_k}) g_k \qfor k\ge 1 .\]
We need to ensure that $\|g_k\|>1$ for all $k \ge 1$.

To this end, let $V_1 = \{ x \in X : h_1(x) > 1 \}$.
Our task will be accomplished if we can find non-empty open sets
$V_{k+1} \subset V_k$ for $k \ge 1$ so that
\[ \sigma_{w_k u_k}(V_{k+1}) \subset V_k \qand
   \sigma_{u_k}(V_k) \subset V_1 \qforal k \ge 1 .
\]

Indeed, if we have these inclusions, we can show by induction that
$g_k > 1$ on $V_k$.
Clearly this is the case for $k=1$.
Assuming the properties outlined in the previous paragraph,
\[ \sigma_{u_k}(V_{k+1}) \subset \sigma_{u_k}(V_k) \subset V_1 ;\]
and thus $h_1\circ \sigma_{u_k} > 1$ on $V_{k+1}$.
Similarly, $\sigma_{w_k u_k}(V_{k+1}) \subset V_k$, and so
$g_k \circ \sigma_{w_k u_k} > 1$ on $V_{k+1}$.
Finally $g_k > 1$ on $V_k \supset V_{k+1}$.
Therefore the product yields $g_{k+1}>1$ on $V_{k+1}$.

To construct $V_{k+1}$, observe that there are no $(t_{k+1},w_1u_k)$-wandering sets.
In particular, there is a $v'_k \in \Fn$ so that
\[ V_{k+1} := \sigma_{t_kv'_kwu_k}^{-1}(V_k) \cap V_k \ne \mt .\]
With this choice (and recalling that $w_{k,1} =t_kv'_kw_1$), we have that
\[ \sigma_{w_{k,1} u_k}(V_{k+1} )\subset V_k \qand V_{k+1} \subset V_k .\]
Lastly,
\[ \sigma_{u_{k+1}}(V_{k+1})
 = \sigma_{u_k} \sigma_{w_{k,1}  u_k}(V_{k+1})
 \subset \sigma_{u_k} (V_k) \subset V_1 .
\]
This completes the induction.
We see that the radical must be $\{0\}$, and so $\A$ is semisimple.

The equivalence of (1) and (4) in the metric case is given by Corollary~\ref{C:metric_wander}.
\end{proof}

\chapter{Open Problems and Future Directions}\label{Problems}

The major goal of associating an operator algebra to a dynamical system
is to provide another venue for seeking invariants for these systems.
In the case of a (amenable) group of homeomorphisms, there have been
a number of successes using C*-algebra crossed products.  
Systems which are not invertible, such as 
those considered here, appear to be better suited to nonself-adjoint
operator algebras.

The major open problem arising from this paper, 
as elaborated in Conjectures \ref{Conj:iso} and \ref{Conj:simplex},
is whether the isomorphism class of the tensor algebra
of a multivariable dynamical systems
is a complete invariant for the system up to piecewise conjugacy.
This is valid for low dimensional topological spaces, and
when the multisystem consists of no more than three maps;
and here isomorphism and completely isometric isomorphism
coincide.
This offers strong evidence for a positive answer.

On the other hand, the classification problem for semicrossed products
is wide open.  We have no plausible conjecture here. 
Example \ref{Ex:pwc but not iso} shows that completely isometric isomorphism
of the semicrossed products can be strictly stronger than
piecewise conjugacy. 
We suspect that, in that example, the algebras are not even algebraically isomorphic.
So what does the isomorphism class of the semicrossed product capture?

Our classification tool, the two dimensional nest representations, 
do not seem to differentiate, in a qualitative way, between the
tensor algebra and the semicrossed product. 
We have no way of extracting from them the additional information 
needed for the classification of the semicrossed products.
Perhaps a study of higher dimensional nest representations, 
as it was done in \cite{DKgraph} for graph algebras,
will yield additional structural information about these operator algebras.

In \cite{KaPet}, it is shown that the fundamental relation for strongly maximal TAF
algebras can be recovered from the (perhaps infinite dimensional)
representation theory of such algebras. 
Since various TAF algebras (as well as our tensor algebras)
are actually tensor algebras of correspondences, we expect 
that the completion of the classification scheme of Muhly and Solel \cite{MS3}
mentioned in the introduction will profit from a better
understanding of their nest representations. 

Also in the case of TAF subalgebras of AF C*-algebras, Power \cite{Pow_lim}
was partially successful in using homology theory to extract invariants.
There may be a role here for some homology groups to capture
topological information. 
\bigbreak

We briefly consider some other aspects of the relationship between the
dynamics and the operator algebras.
In Chapter~\ref{C:semisimple}, we clarified the connection between
the semisimplicity of the operator algebra and the (non-existence of)
generalized wandering sets.
This naturally begs the question of a topological description of the
Jacobson radical in terms of the underlying dynamics.
In the case of a single map, 
or even a finitely generated family of commuting maps, 
this was done in \cite{DKM}. 
Likely, as in their work, a refinement of the notion of recurrence will
play a role.

One can also consider other topological conditions such as chaos,
and try to capture this in the operator algebra context.
In the case of a single map $\sigma$, one can observe the
following phenonena.
A periodic point leads to a finite dimensional \textit{maximal} 
orbit representation.  So a dense set of periodic points means that
the C*-envelope has a faithful representation which is the 
direct sum of irreducible finite dimensional representations.
In particular, it is quasidiagonal.
What does quasidiagonality of this C*-algebra mean in general?
Also a transitive point yields a faithful irreducible maximal representation;
that is, the C*-envelope is primitive.
Do these two properties together imply that the system is chaotic?
Is there an analogue for multivariable systems?
\bigbreak

A third aspect that we wish to mention is consideration of relations
among the maps $\sigma_i$.  It is possible to encode certain types
of relations and create a universal operator algebra.  Such relations
will be considered in a forthcoming paper.  
These relations have the effect of restricting the possible characters
to lie on an analytic variety in $\bC^n$.  So it is anticipated that
further connections to several complex variables or algebraic
geometry will be required.
One is also tempted to consider local relations,
such as maps commuting when the domain lies in a certain open subset. 
A much more elusive quest is whether one can construct,
in a natural way, an operator algebra that automatically
encodes any such relationships that exist.

\backmatter

\end{document}